\documentclass[11pt,leqno]{article}

\usepackage[dvips]{graphicx}

\usepackage{amsmath}
\usepackage{amsthm}
\usepackage{amsfonts}

\reversemarginpar
\newcommand{\mmark}[1]{\marginpar{\footnotesize\tt{#1}}}
\newcommand{\XXX}[1]{\marginpar{\footnotesize\tt{#1}}}

\renewcommand{\mmark}[1]{}
\renewcommand{\XXX}[1]{}

\newcommand{\commentout}[1]{}

\theoremstyle{plain}
\newtheorem{theorem}{Theorem}[section]
\newtheorem{proposition}[theorem]{Proposition}
\newtheorem{lemma}[theorem]{Lemma}
\newtheorem{corollary}[theorem]{Corollary}
\newtheorem{conjecture}[theorem]{Conjecture}

\theoremstyle{definition}
\newtheorem{definition}[theorem]{\it Definition}
\newtheorem{example}[theorem]{\it Example}
\newtheorem{statement}[theorem]{\it Statement}

\theoremstyle{remark}
\newtheorem{remark}[theorem]{Remark}

\newcommand{\SSection}[2][??]{
            \section{#2}
            \label{SS:#1}
            \mmark{SS:#1}
            \setcounter{equation}{0}
            \setcounter{table}{0}
            \setcounter{figure}{0}}
\newcommand{\SubSSection}[2][??]{
            \subsection{#2}
            \label{SS:#1}
            \mmark{SS:#1}}

\newcommand{\SStop}{\hspace\parindent}

\newcommand{\References}{\newpage\begin{thebibliography}{99}}
\newcommand{\EndReferences}{\end{thebibliography}}

\newcommand{\SSn}[1]{\ref{SS:#1}}

\newcommand{\Proposition}[1][??]{
            \begin{proposition}
            \label{TH:#1}
            \mmark{TH:#1}}
\newcommand{\EndProposition}{\end{proposition}}

\newcommand{\Theorem}[1][??]{
            \begin{theorem}
            \label{TH:#1}
            \mmark{TH:#1}}
\newcommand{\EndTheorem}{\end{theorem}}

\newcommand{\Lemma}[1][??]{
            \begin{lemma}
            \label{TH:#1}
            \mmark{TH:#1}}
\newcommand{\EndLemma}{\end{lemma}}

\newcommand{\Corollary}[1][??]{
            \begin{corollary}
            \label{TH:#1}
            \mmark{TH:#1}}
\newcommand{\EndCorollary}{\end{corollary}}

\newcommand{\Definition}[1][??]{
            \begin{definition}
            \label{TH:#1}
            \mmark{TH:#1}}
\newcommand{\EndDefinition}{\end{definition}}

\newcommand{\Conjecture}[1][??]{
            \begin{conjecture}
            \label{TH:#1}
            \mmark{TH:#1}}
\newcommand{\EndConjecture}{\end{conjecture}}

\newcommand{\Example}[1][??]{
            \begin{example}
            \label{TH:#1}
            \mmark{TH:#1}}
\newcommand{\EndExample}{\end{example}}

\newcommand{\Remark}[1][??]{
            \begin{remark}
            \label{TH:#1}
            \mmark{TH:#1}}
\newcommand{\EndRemark}{\end{remark}}

\newcommand{\Statement}[1][??]{
            \begin{statement}
            \label{TH:#1}
            \mmark{TH:#1}}
\newcommand{\EndStatement}{\end{statement}}

\newcommand{\THn}[1]{\ref{TH:#1}}

\newcommand{\Proof}[1][{\bf Proof.}]{\begin{proof}[#1]}
\newcommand{\QED}{\end{proof}}

\renewcommand{\qed}{\setlength{\unitlength}{1.0ex}
  \begin{picture}(1.4,1.4)
    \put(0,0){\framebox(1.4,1.4)}
  \end{picture}}

\newcommand{\Equation}[1][??]{
            \mmark{EQ:#1}
            \begin{equation}
            \label{EQ:#1}}

\newcommand{\EQn}[1]{{\rm (\ref{EQ:#1})}}

\newcounter{abcItem}
\newcommand{\abcList}{\begin{list}{{\rm (\alph{abcItem})}}{
  \usecounter{abcItem}
  \setlength{\parsep}{0.2ex}
}}
\newcounter{numItem}
\newcommand{\numList}{\begin{list}{{\rm (\arabic{numItem})}}{
  \usecounter{numItem}
  \setlength{\parsep}{0.2ex}
}}
\newcounter{romanItem}
\newcommand{\romanList}{\begin{list}{{\rm (\roman{romanItem})}}{
  \usecounter{romanItem}
  \setlength{\parsep}{0.2ex}
}}
\newcommand{\anyList}[1]{\begin{list}{#1}{
  \setlength{\parsep}{0.2ex}
}}
\newcommand{\contList}{\begin{list}{}{
  \setlength{\parsep}{0.2ex}
  \setlength{\topsep}{0.0ex}
  \setlength{\labelwidth}{1.0ex}
  \setlength{\labelsep}{0.0ex}
  \setlength{\itemindent}{0.0ex}
}}
\newcommand{\EndList}{\end{list}}

\newcommand{\Ac}{{\mathcal A}}
\newcommand{\Bc}{{\mathcal B}}
\newcommand{\Cc}{{\mathcal C}}

\newcommand{\Ec}{{\mathcal E}}
\newcommand{\Fc}{{\mathcal F}}
\newcommand{\Gc}{{\mathcal G}}

\newcommand{\Mc}{{\mathcal M}}

\newcommand{\Pc}{{\mathcal P}}

\newcommand{\Sc}{{\mathcal S}}

\newcommand{\Uc}{{\mathcal U}}
\newcommand{\Vc}{{\mathcal V}}
\newcommand{\Wc}{{\mathcal W}}

\newcommand{\PP}{{\bf P}}

\newcommand{\UU}{{\bf U}}
\newcommand{\XX}{{\bf X}}

\newcommand{\xx}{{\bf x}}

\newcommand{\B}{{\rm B}}

\newcommand{\W}{{\rm W}}
\newcommand{\Y}{{\rm Y}}
\newcommand{\Z}{{\rm Z}}

\newcommand{\Prob}{\mathbb P}
\newcommand{\stocle}{\preceq}

\newcommand{\unif}{\textrm{unif}}

\newcommand{\pideal}[1]{\langle{#1}\rangle}

\newcommand{\csection}[1]{(\leftarrow,{#1}]}

\begin{document}

\title{
\commentout{
\fbox{\parbox[b]{5.3in}{
\begin{center} {\huge\textsf{DRAFT}} \end{center}
{\normalsize \begin{itemize}
\item \textsf{DRAFT} may contain inaccurate expositions to undergo
further editing
\item Revisions will be posted at
{\tt http://www.mts.jhu.edu/\symbol{126}fill/}
\end{itemize}}
}}
\\[0.3in]}
Stochastic Monotonicity and\\
Realizable Monotonicity
}

\author{
     James Allen Fill
\thanks{Research for both authors was supported by
        NSF grants DMS-96-26756
        and DMS-98-03780.}
\and Motoya Machida
\and \\
     Department of Mathematical Sciences\\
     The Johns Hopkins University
}

\date{\today}

\maketitle

\begin{abstract}
We explore and relate two notions of monotonicity, stochastic and
realizable, for a system of probability measures on a common finite
partially ordered set (poset) $\Sc$ when the measures are indexed by
another poset $\Ac$.
We give counterexamples to show that the two notions are not
always equivalent, but for various large classes of $\Sc$ we
also present conditions on the poset $\Ac$
that are necessary and sufficient for equivalence.
When $\Ac = \Sc$, the condition that
the cover graph of $\Sc$ have no cycles is necessary and sufficient
for equivalence.
This case arises in comparing applicability of the perfect sampling
algorithms of Propp and Wilson and the first author of the present paper.

\vspace{1.5in}

\par\noindent
{\em Short title.\/}
Stochastic and realizable monotonicity.

\medskip
\par\noindent
{\em AMS\/} 1991 {\em subject classifications.\/}
Primary 60E05;
secondary 06A06, 60J10, 05C38.


\end{abstract}

\pagebreak

\SSection[intro]{Introduction}

\SStop
We will discuss two notions of monotonicity for probability measures on 
a finite partially ordered set (poset).
Let $\Sc$ be a finite poset and
let $(P_1,P_2)$ be a pair of probability measures on $S$.
(We use a calligraphic letter $\Sc$ in order to distinguish the set
$S$ from the same set equipped with a partial ordering $\le$.)
A subset $U$ of $S$ is said to be an {\em up-set\/} in $\Sc$
(or {\em increasing set\/})
if $y \in U$ whenever $x \in U$ and $x \le y$.
We say that $P_1$ is {\em stochastically smaller\/} than $P_2$,
denoted
\Equation[stocle]
  P_1 \stocle P_2,
\end{equation}
if
\Equation[stocle.def]
P_1(U) \le P_2(U)
\quad\text{ for every up-set $U$ in $\Sc$.}
\end{equation}
The relation introduced in~\EQn{stocle}--\EQn{stocle.def} is clearly reflexive and transitive.
Antisymmetry follows easily using our assumption that $S$ is finite,
so the relation defines a partial ordering on the class of probability
measures on $S$.
(For a careful discussion on the matter of antisymmetry in a rather
general setting for infinite $S$, see~\cite{KK}.)

The following characterization of stochastic ordering was established
by Strassen~\cite{Strassen} and fully investigated by 
Kamae, Krengel, and O'Brien~\cite{KKO}.
Suppose that there exists a pair $(\XX_1, \XX_2)$ of $S$-valued random variables
[defined on some probability space $(\Omega,\Fc,\Prob)$]
satisfying the properties
\Equation[rv.le]
  \XX_1 \le \XX_2
\end{equation}
and
\Equation[rv.marg]
  \Prob(\XX_i \in \cdot) = P_i(\cdot)
  \quad\text{ for $i=1,2$. }
\end{equation}
Then we have
$$
P_1(U) = \Prob(\XX_1\in U) = \Prob(\XX_1\in U,\:\XX_1\le\XX_2)
\le \Prob(\XX_2\in U) = P_2(U),
$$
for every up-set $U$ in $\Sc$.
Thus, the conditions \EQn{rv.le}--\EQn{rv.marg} necessitate
\EQn{stocle}.
Moreover, Strassen's work shows that \EQn{stocle} is in fact sufficient for the existence of
a probability space $(\Omega,\Fc,\Prob)$ and
a pair $(\XX_1, \XX_2)$ of $S$-valued random variables on $(\Omega,\Fc,\Prob)$
satisfying \EQn{rv.le}--\EQn{rv.marg}.
[Equivalently, we need only require that \EQn{rv.le} hold almost surely.]

Now let $\Ac$ be a finite poset.
Let $(P_\alpha: \alpha\in A)$ be a system of
probability measures on $S$.
We call $(P_\alpha: \alpha\in A)$ a {\em realizably monotone\/} system
if there exists a system $(\XX_\alpha: \alpha\in A)$
of $S$-valued random variables defined on some probability space
$(\Omega,\Fc,\Prob)$ such that
\Equation[rm.mono]
  \XX_\alpha \le \XX_\beta
  \quad\text{ whenever $\alpha \le \beta$ }
\end{equation}
and
\Equation[rm.marg]
\Prob(\XX_\alpha \in \cdot) = P_\alpha(\cdot)
\quad\text{ for every $\alpha\in A$. }
\end{equation}
In such a case we shall say that $(\XX_\alpha: \alpha\in A)$
{\em realizes the monotonicity\/} of $(P_\alpha: \alpha\in A)$.
Since the conditions~\EQn{rv.le}--\EQn{rv.marg} imply~\EQn{stocle},
the conditions~\EQn{rm.mono}--\EQn{rm.marg} applied pairwise imply
\Equation[sm.mono]
  P_\alpha \stocle P_\beta \quad\text{ whenever $\alpha\le\beta$. }
\end{equation}
The system $(P_\alpha: \alpha\in A)$ is said to be {\em stochastically monotone\/}
if it satisfies~\EQn{sm.mono}.
We have shown that stochastic monotonicity is necessary for realizable
monotonicity.

In light of Strassen's characterization of stochastic ordering,
one might guess that stochastic monotonicity is also sufficient for realizable monotonicity.
It is perhaps surprising that the conjecture is false in general,
as the following example shows.

\Example[ex.dd]
Let
\Equation[diamond]
\Sc  = \Ac :=
\begin{minipage}{0.8in}\begin{center}
  \includegraphics{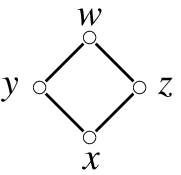}
\end{center}\end{minipage}
\end{equation}
be the usual 2-dimensional Boolean algebra
with $x < y, z$ and $y, z < w$
(and, of course, $x < w$ by transitivity).
Define a system $(P_x,P_y,P_z,P_w)$ of probability measures on $S$ by
$$
  P_\xi := \begin{cases}
    \unif\{x,y\} & \quad\text{ if $\xi = x$; } \\
    \unif\{x,w\} & \quad\text{ if $\xi = y$; } \\
    \unif\{y,z\} & \quad\text{ if $\xi = z$; } \\
    \unif\{y,w\} & \quad\text{ if $\xi = w$, }
  \end{cases}
$$
where $\unif(B)$ denotes the uniform probability measure on a set $B$.
Clearly, $(P_x,P_y,P_z,P_w)$ is stochastically monotone.
Now suppose that there exists a system $(\XX_x,\XX_y,\XX_z,\XX_w)$
which realizes the monotonicity of $(P_x,P_y,P_z,P_w)$.
Considering the event $\XX_x = y$, realizable monotonicity forces
$$
  \Prob(\XX_x = y) = \Prob(\XX_{x} = y,\:\XX_{y} = w,\:\XX_{z} = y,\:\XX_{w} = w) = \frac{1}{2}.
$$
Similarly, we find
$$
  \Prob(\XX_z = z) = \Prob(\XX_{x} = x,\:\XX_{z} = z,\:\XX_{w} = w) = \frac{1}{2}.
$$
Noting that the above two events are disjoint,
we conclude $\Prob(\XX_{w}=w) = 1$, which is a contradiction.
Thus, $(P_x,P_y,P_z,P_w)$ cannot be realizably monotone.
\qed
\EndExample

Given a pair $(\Ac,\Sc)$ of posets,
if the two notions of monotonicity---stochastic and realizable---are
equivalent, then we say that {\em monotonicity equivalence\/} holds
for $(\Ac,\Sc)$.
The counterexample in Example~\THn{ex.dd} was discovered
independently by Ross~\cite{Ross};
we are grateful to Robin Pemantle for pointing this out to us.
Ross reduced the question of monotonicity equivalence for general
{\em infinite\/} posets $\Ac$ (and given $\Sc$) to consideration of
the same question for every finite induced subposet of $\Ac$.
Thus we regard our work as a useful complement to his.
As an historical aside, we note that what we call a realizably
monotone system, Ross called a {\em coherent family\/}.

To give an example where monotonicity equivalence holds, we next
consider the case where $\Sc$ is a {\em linearly\/} ordered set.

\Example[ex.linear]
Let $\Ac$ be any poset and let $\Sc$ be a linearly ordered set.
Suppose that $(P_\alpha: \alpha\in A)$ is a stochastically monotone system of
probability measures on $S$.
For each $\alpha\in A$, define the {\em inverse probability transform\/}
$P_\alpha^{-1}$ by
\Equation[def.inv.pr]
  P_\alpha^{-1}(t)
  := \min\left\{x\in S: t < F_\alpha(x) \right\}
  \quad\text{ for $t \in [0,1)$,}
\end{equation}
where $F_\alpha$ is the distribution function of $P_\alpha$
[i.e., $F_\alpha(x) := P_\alpha(\{\xi\in S: \xi\le x\})$ for each $x \in S$].
Given a single uniform random variable $\UU$ on $[0,1)$,
we can construct a system $(\XX_\alpha: \alpha\in A)$ of
$S$-valued random variables via
$$
  \XX_\alpha := P_\alpha^{-1}(\UU)
  \quad\text{ for each $\alpha\in A$.}
$$
Then $(\XX_\alpha: \alpha\in A)$ realizes the monotonicity
and therefore $(P_\alpha: \alpha\in A)$ is realizably monotone.
Thus monotonicity equivalence holds for $(\Ac,\Sc)$.
\qed
\EndExample

The goal of our investigation is to determine for
precisely which pairs $(\Ac,\Sc)$ of posets monotonicity equivalence
holds.
Let us discuss here the usefulness of such a determination.
It is (structurally) simple to say which systems
$(P_\alpha: \alpha\in A)$ are stochastically monotone.
Indeed, one need only determine all up-sets $U$ of $\Sc$,
and then $(P_\alpha: \alpha\in A)$ is stochastically monotone
if and only if $P_\alpha(U) \le P_\beta(U)$ for all such $U$
whenever $\alpha\le\beta$.
For realizable monotonicity, an analogous result is
Theorem~\THn{ext.thm}, but the necessary and
sufficient condition there involves an infinite collection of
inequalities.
We know how to reduce, for each $(\Ac,\Sc)$, the infinite collection
to a finite one, but (1)~there seems to be in general no nice
structural characterization of the resulting finite collection,
and (2) the computations needed to do the reduction can be massive
even for fairly small $\Ac$ and $\Sc$
(Chapter~7 of~\cite{mmthesis}).

Thus, when monotonicity equivalence holds, we learn that realizable
monotonicity has the same simple structure as stochastic
monotonicity.
And when monotonicity equivalence fails, we learn that testing a
system $(P_\alpha: \alpha\in A)$ for stochastic monotonicity does
{\em not\/} suffice as a test for realizable monotonicity.
In the latter case, for fixed $(\Ac,\Sc)$ and a single numerically
specified system $(P_\alpha: \alpha\in A)$, we can determine whether
or not the system is realizably monotone by constructing a system
$(\XX_\alpha: \alpha\in A)$ subject to the marginal
condition~\EQn{rm.marg} so as to maximize the probability
$\Prob(\XX_\alpha \le \XX_\beta \text{ whenever $\alpha\le\beta$})$.
Indeed, the system $(P_\alpha: \alpha\in A)$ is realizably monotone if 
and only if the maximum value equals $1$.
The construction can be carried out using linear programming with
variables corresponding to the values of the joint probability mass
function for $(\XX_\alpha: \alpha\in A)$.

For further discussion along these lines, see~\cite{mmthesis}.

\begin{figure}[h]
\vspace{0.2in}
\begin{center}
\begin{tabular}{ccc}
$\Sc :=$
\begin{minipage}{2.3in}\begin{center}
  \includegraphics{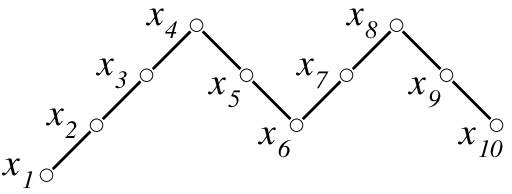}
\end{center}\end{minipage}
& \hspace{0.2in} &
$(S,\le_{n}) :=$
\begin{minipage}{0.7in}\begin{center}
  \includegraphics{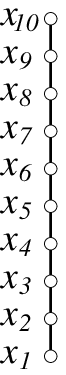}
\end{center}\end{minipage} \\[0.7in]
(a) a poset of Class~Z & & (b) a linear order of the path
\end{tabular}
\caption{The comparison of two posets}
\label{fig.class.z}
\end{center}
\end{figure}

Of particular interest in our study of realizable monotonicity
is the case $\Ac = \Sc$.
Here the system $(P(x,\cdot): x \in S)$
of probability measures can be considered as a Markov
transition matrix $\PP$ on the state space $S$.
Recently, Propp and Wilson~\cite{Propp-Wilson} and Fill~\cite{Fill}
have introduced algorithms to produce observations distributed
{\em perfectly\/} according to the long-run distribution of a Markov
chain.
Both algorithms apply most readily and operate most efficiently when
the state space $\Sc$ is a poset and a suitable monotonicity condition
holds.
Of the many differences between the two algorithms, one is that the
appropriate notion of monotonicity for the Propp--Wilson algorithm is
realizable monotonicity, while for Fill's algorithm it is stochastic
monotonicity; see Remark~4.5 in~\cite{Fill}.
Here the properties~\EQn{rm.mono}--\EQn{rm.marg} are essential
for the Propp--Wilson algorithm to be able to generate
transitions simultaneously from every state in such a way
as to preserve ordering relations.
For further discussion of these perfect sampling algorithms,
see~\cite{Fill} and~\cite{Propp-Wilson}.
In Theorem~\THn{me.markov} we show that the two notions of
monotonicity are equivalent if and only if the poset $\Sc$ is
``acyclic'',
which is characterized by
possession of a Hasse diagram
(the standard graphical representation of partial ordering)
that is cycle-free.
For example, the Hasse diagram of a linearly ordered set is a
vertical path such as the one in Figure~\ref{fig.class.z}(b), and therefore
has no cycle.
On the other hand, the 2-dimensional Boolean algebra whose Hasse diagram is
displayed in~\EQn{diamond} is not acyclic.
See Section~\SSn{poset} for precise terminology.

In the present paper we study the notion of realizable
monotonicity when $\Ac$ and $\Sc$ are both finite posets.
In Section~\SSn{strassen}
we review a general result for the existence of a probability measure
with specified marginals and present the extensibility problem.
Sections~\SSn{poset}--\SSn{cycle} are prepared to introduce
definitions and several key notions in studying posets.
In Section~\SSn{me.problem}
we formulate the monotonicity equivalence problem
from the viewpoint of the extensibility problem.
Section~\SSn{me.summary} is rather short,
introducing four subclasses---Classes~B,~Y,~W~and~Z---that
partition the class of connected posets $\Sc$.

At this juncture we intend to provide the reader with an
overview of the main results of this paper.
In Section~\SSn{class.b} we present the first case of our investigation,
where $\Sc$ is in Class~B.
Kamae,\ Krengel,\ and O'Brien~\cite{KKO} showed that if $\Ac$ is a
linearly ordered set then monotonicity equivalence holds for
$(\Ac,\Sc)$.
We generalize this result (in our finite setting) to the case of an acyclic
poset $\Ac$ (Theorem~\THn{a.acyclic});
see Section~\SSn{poset} for the definition of an acyclic poset.
Theorem~\THn{class.b.main} gives an exact answer to our central question
of monotonicity equivalence when $\Sc$ is a poset of Class~B.
In Section~\SSn{class.y} we proceed to the second case of our investigation,
where $\Sc$ is in Class~Y.
There we first present Proposition~\THn{p.insert} and discuss its proof in
Section~\SSn{prob.on.acyc}.
This turns out to be a useful result concerning probability measures on a poset of
Class~Y, leading to Theorem~\THn{class.y.main},
which in turn answers our monotonicity equivalence question when $\Sc$ is a poset of Class~Y.
If $\Sc$ is a poset of Class~Z,
then we can show that monotonicity equivalence holds for any poset $\Ac$
(Theorem~\THn{class.z.main}).
In Section~\SSn{class.z} we give the proof by using a generalization of inverse
probability transform.
When $\Sc$ is a poset of Class~W,
we have devised a further generalization of inverse probability transform,
which results in constructing a rather large class of posets
$\Ac$ for which monotonicity equivalence holds.
We refer the reader to~\cite{invprob} for the results of our investigation
of Class~W.

\SSection[ext.prob]{Posets and the monotonicity equivalence problem}

\SStop
In Section~\SSn{poset} we briefly summarize the material on
posets that we need for our study.
An important assertion is that if a
poset is non-acyclic then the poset has an induced cyclic subposet.
In Section~\SSn{cycle}
we prove this (Lemma~\THn{cons.cycle}) among other results concerning
induced cyclic posets.
In Section~\SSn{strassen} we review the well-known results
of Strassen~\cite{Strassen} on the
existence of a probability measure with specified marginals;
our review is tailored somewhat to fit our application to realizable
monotonicity.
In Section~\SSn{me.problem}
we discuss realizable monotonicity
in terms of the existence of a probability measure with specified marginals.
Propositions~\THn{a.connect}--\THn{s.connect}
are presented in Section~\SSn{me.problem};
these allow $\Ac$ and $\Sc$ both to be
connected posets in our later investigations.

\SubSSection[poset]{Posets}

\SStop
We devote this subsection to introducing definitions and
notation related to partial ordering.
By a {\em poset\/} $\Sc$ we shall mean a {\em finite\/} set $S$
(the qualifier will not again be explicitly noted)
together with a partial ordering $\le$.
The (unordered) set $S$ is called the {\em ground set\/} of $\Sc$.
Most of the basic poset terminology adopted here can be found
in Stanley~\cite{Stanley} or Trotter~\cite{Trotter}.
Throughout this subsection, $\Sc$ and $\Sc'$ denote posets.

\numList
\item
{\em Dual poset, up-set, down-set.\/}
The {\em dual} of $\Sc$, denoted $\Sc^*$,
is the poset on the same ground set $S$ as $\Sc$ such that
$x_1 \le x_2$ in $\Sc^*$ if and only if $x_1 \ge x_2$ in $\Sc$.
A subset $U$ of $S$ is said to be an {\em up-set\/}
(or {\em increasing set\/}) in $\Sc$
if $y \in U$ whenever $x \in U$ and $x \le y$.
A {\em down-set} $V$ in $\Sc$ is defined to be an up-set in $\Sc^*$.
Note that $U$ is an up-set in $\Sc$
if and only if $S \setminus U$ is a down-set in $\Sc$.
For any subset $B$ of $S$, we can define the {\em down-set $\pideal{B}$
generated by\/} $B$ in the usual fashion:
$$
  \pideal{B}
  := \{\xi\in S: \xi \le \eta \text{ for some $\eta\in B$}\}.
$$
We simply write $\pideal{x_1,\ldots,x_n}$ for $\pideal{\{x_1,\ldots,x_n\}}$.

\item
{\em Cover graph.\/}
For $x, y \in S$,
we say that $y$ {\em covers\/} $x$ if $x < y$ in $\Sc$ and no element
$z$ of $S$ satisfies $x < z < y$.
The {\em Hasse diagram\/} of a poset is the directed graph whose vertices
are the elements of the poset and whose arcs are those ordered pairs $(x,y)$
such that $y$ covers $x$.
[By convention, if $y$ covers $x$, then $y$ is drawn above $x$
in the Hasse diagram (as represented in the plane);
this indicates the direction of each arc.]
We define the {\em cover graph\/} $(S,\Ec_{\Sc})$ of $\Sc$ by considering the
Hasse diagram of $\Sc$ as an undirected graph.
That is, the edge set $\Ec_{\Sc}$ consists of those unordered pairs $\{x,y\}$
such that either $x$ covers $y$ or $y$ covers $x$ in $\Sc$.

\item
{\em Subposet.\/}
We shall need to distinguish among several, somewhat subtly different,
notions of subposet.
We say that a poset $\Sc'$ is
a {\em subposet\/} of $\Sc$
if $S'$ is
(or, by extension and when there is no possibility of confusion, is
isomorphic to)
a subset of $S$ and
$x \le y$ in $\Sc'$ implies $x \le y$ in $\Sc$ for $x, y \in S'$.
When we speak of an {\em induced subposet\/} $\Sc'$ of $\Sc$,
we mean that for $x, y \in S'$ we have $x \le y$ in $\Sc'$
if and only if $x \le y$ in $\Sc$.
On the other hand, we call a (not necessarily induced) subposet $\Sc'$
a {\em subposet via induced cover subgraph\/} of $\Sc$
if $y$ covers $x$ in $\Sc'$ for $x, y \in S'$ precisely when $y$
covers $x$ in $\Sc$, that is,
when the cover graph $(S',\Ec_{\Sc'})$ of $\Sc'$ is an induced subgraph of
the cover graph $(S,\Ec_{\Sc})$ of $\Sc$.
Clearly, a subposet via induced cover subgraph of $\Sc$ with ground set $S'$
is a subposet of the subposet of $\Sc$ induced by ground set $S'$.
In Example~\THn{ex.induce} we illustrate differences between these two
notions of subposet.
\EndList

Let $(S',\Ec')$ be a (not necessarily induced) subgraph of
the cover graph $(S,\Ec_{\Sc})$ of $\Sc$.
Then $(S',\Ec')$ is the cover graph $(S',\Ec_{\Sc'})$
of a (not necessarily induced) subposet $\Sc'$ of $\Sc$.
Here, $y$ covers $x$ in $\Sc'$ if and only if $y$ covers $x$ in $\Sc$ and
$\{x,y\}\in\Ec'$.
In this sense, a subgraph $(S',\Ec')$ of $(S,\Ec_{\Sc})$
can be considered as a subposet of $\Sc$.

\numList
\setcounter{numItem}{3}
\item
{\em Chain, height.\/}
We call a poset $\Sc$ a {\em chain} if any two elements of $S$ are
comparable in $\Sc$.
When we say that a subposet $\Sc'$ is a chain in $\Sc$,
we mean that $\Sc'$ is a chain and an induced subposet of $\Sc$.
The {\em height\/} $n$ of a poset $\Sc$ is the number of elements in
a maximum-sized chain in $\Sc$.
That is, $\Sc$ has height $n$ if and only if
$\Sc$ has an $n$-element chain, but no $(n+1)$-element chain,
as an induced subposet.

\item
{\em Path, upward path, downward path.\/}
We call a (not necessarily induced) subposet $\Sc'$ of $\Sc$
a {\em path\/} if the cover graph $(S',\Ec_{\Sc'})$ of $\Sc'$ is
(i) a path (in the usual graph-theoretic sense)
and (ii) a (not necessarily induced) subgraph of the cover graph
$(S,\Ec_{\Sc})$ of $\Sc$.
A sequence $(x_0,x_1,\ldots,x_{n-1})$ denotes a path from $x_0$ to $x_{n-1}$
with vertex set $\{x_0,x_1,\ldots,x_{n-1}\}$ and
edge set $\{\{x_{i-1},x_{i}\}: 1 \le i \le n-1\}$.
We say that a path $(x_0,\ldots,x_{n-1})$ is {\em upward\/}
(respectively, {\em downward\/}) in $\Sc$ if
$x_{i}$ covers $x_{i-1}$ in $\Sc$ ($x_{i-1}$ covers $x_{i}$, respectively)
for each $i=1,\ldots,n-1$.
Note that any upward or downward path in $\Sc$ is a chain in $\Sc$,
but that the converse is not true.
We illustrate chains and paths in Example~\THn{ex.induce}(iv)--(vi).

\item
{\em Cycle.\/}
We call a (not necessarily induced) subposet $\Sc'$ of $\Sc$
a {\em cycle\/} (or a {\em cyclic subposet\/})
if the cover graph $(S',\Ec_{\Sc'})$ of $\Sc'$ is
(i) a cycle (in the usual graph-theoretic sense)
and (ii) a (not necessarily induced) subgraph of the cover graph
$(S,\Ec_{\Sc})$ of $\Sc$.
In Example~\THn{ex.induce} we demonstrate that a cyclic
subposet may be neither an induced subposet nor
a subposet via induced cover subgraph.
If a poset $\Sc$ has a cyclic subposet, then
we call $\Sc$ a {\em non-acyclic\/} poset.
In keeping with the foregoing definitions, we call the reference poset
$\Sc$ itself a {\em cycle\/} if the cover graph $(S,\Ec_{\Sc})$ of
$\Sc$ is a cycle.
A sequence $(x_0,x_1,\ldots,x_{n-1},x_0)$ with $n \ge 4$
denotes a cycle with vertex set $\{x_0,x_1,\ldots,x_{n-1}\}$ and
edge set $\{\{x_{i-1},x_{i}\}: 1 \le i \le n\}$.
Here, indices are interpreted modulo $n$.
(Note that a cyclic subposet must consist of at least four elements.)

\item
{\em Connected poset, disjoint union.\/}
We say that $\Sc$ is {\em connected\/} if its cover graph $(S,\Ec_{\Sc})$
is connected.
The {\em components\/} of $\Sc$ are its maximal connected induced
subposets.
If $S$ and $S'$ are disjoint, then
we can construct the {\em disjoint union\/} of $\Sc$ and $\Sc'$,
denoted $\Sc + \Sc'$,
as a poset on the ground set $S \cup S'$ by declaring $x \le y$ in $\Sc + \Sc'$
precisely when either (i) $x,y \in S$ and $x \le y$ in $\Sc$,
or (ii) $x,y \in S'$ and $x \le y$ in $\Sc'$.
Thus any poset $\Sc$ is the disjoint union of its components.

\item
{\em Acyclic poset, leaf.\/}
We say that a poset $\Sc$ is {\em acyclic\/} if $\Sc$ has no cyclic subposet.
We call an element $x$ of $S$ a {\em leaf\/} in $\Sc$ if
the edge set $\Ec_{\Sc}$ of the cover graph of $\Sc$
has a unique element $\{x,y\}$ for some $y \in S$.
Note that if $x$ is a leaf in $\Sc$ then $x$ must be either maximal or
minimal in $\Sc$.
If $\Sc$ is a connected acyclic poset with $|S| \ge 2$, then
there are at least two leaves in $\Sc$ (see, e.g.,~\cite{West}).

\begin{figure}[h]
\vspace{0.2in}
\begin{center}
\begin{tabular}{ccccc}
\begin{minipage}{0.8in}\begin{center}
  \includegraphics{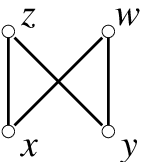}
\end{center}\end{minipage}
& \hspace{0.2in} &
\begin{minipage}{0.8in}\begin{center}
  \includegraphics{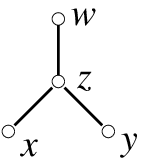}
\end{center}\end{minipage}
& \hspace{0.2in} &
\begin{minipage}{0.8in}\begin{center}
  \includegraphics{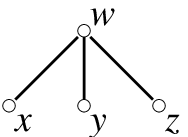}
\end{center}\end{minipage} \\[0.3in]
(a) Bowtie & & (b) a Y-poset & & (c) a W-poset
\end{tabular}
\\[0.3in]
\begin{minipage}{2.4in}\begin{center}
\includegraphics{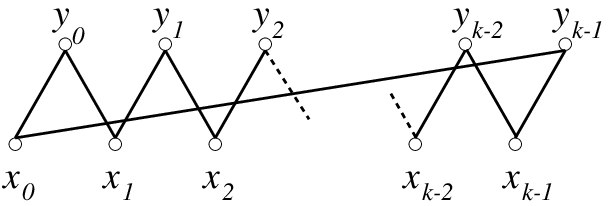}
\end{center}\end{minipage} \\[0.1in]
(d) $k$-crown
\caption{Some named posets}
\label{fig.posets}
\end{center}
\end{figure}

\item
{\em Poset-isomorphism, some named posets, subdivision.\/}
$\Sc$ is said to be {\em poset-isomorphic\/} to $\Sc'$
if there exists an bijection $\phi$ from $S$ to $S'$ such that
$x \le y$ in $\Sc$ if and only if $\phi(x) \le \phi(y)$ in $\Sc'$.
In this paper, we call the $2$-dimensional Boolean algebra
a {\em diamond\/}.
[See the figure in~\EQn{diamond}.]
Furthermore, we call the posets of Figure~\ref{fig.posets} and their duals
(a) the {\em bowtie}, (b) {\em Y-posets}, (c) {\em W-posets},
and (d) the {\em $k$-crown\/}, respectively.
The bowtie is the same as the $2$-crown.
We may simply call $\Sc$ a {\em crown\/} if $\Sc$ is the $k$-crown for
some $k \ge 2$.
We say that a poset $\Sc'$ is a {\em subdivision\/} of $\Sc$
if the induced subposet of $\Sc'$ on ground set
$\{z\in S': x \le z \le y \text{ in $\Sc'$}\}$ is a chain
whenever $y$ covers $x$ in $\Sc$.

\EndList

\Example[ex.induce]
Let
\begin{center}
$\Sc := $
\begin{minipage}{1.0in}\begin{center}
  \includegraphics{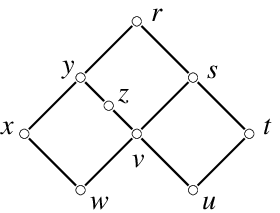}
\end{center}\end{minipage}
\end{center}
be a poset.
Here we give a number of examples to illustrate subtle distinctions
in the definitions of subposets, paths, and cycles.
Let
\begin{center}
$\Sc'_1 := $
\begin{minipage}{1.0in}\begin{center}
  \includegraphics{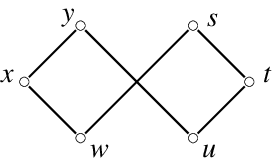}
\end{center}\end{minipage}
\hspace{0.2in}
$\Sc'_2 := $
\begin{minipage}{1.0in}\begin{center}
  \includegraphics{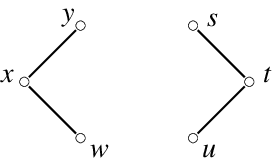}
\end{center}\end{minipage}
\hspace{0.1in} and \hspace{0.1in}
$\Sc'_3 := $
\begin{minipage}{1.0in}\begin{center}
  \includegraphics{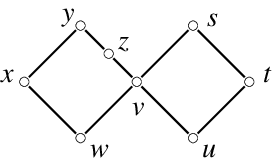}
\end{center}\end{minipage}
\end{center}
be subposets of $\Sc$.
Then, (i) $\Sc'_1$ is an induced subposet but not a subposet via
induced cover subgraph,
(ii) $\Sc'_2$ is a subposet via induced cover subgraph but not an
induced subposet, and
(iii) $\Sc'_3$ is both an induced subposet and
a subposet via induced cover subgraph.
Let
\begin{center}
$\Sc'_4 := $
\begin{minipage}{1.0in}\begin{center}
  \includegraphics{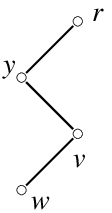}
\end{center}\end{minipage}
\hspace{0.2in}
$\Sc'_5 := $
\begin{minipage}{1.0in}\begin{center}
  \includegraphics{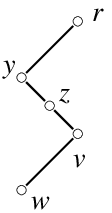}
\end{center}\end{minipage}
\hspace{0.1in} and \hspace{0.1in}
$\Sc'_6 := $
\begin{minipage}{1.0in}\begin{center}
  \includegraphics{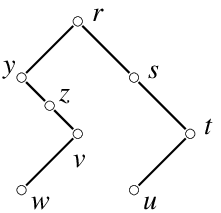}
\end{center}\end{minipage}
\end{center}
be subposets of $\Sc$.
Then, (iv) $\Sc'_4$ is a chain but not a path in $\Sc$,
(v) $\Sc'_5$ is a chain and an upward path from $w$ to $r$,
and (vi) $\Sc'_6$ is a path between $w$ and $u$
but neither an upward path nor a downward path.
Let
\begin{center}
$\Sc'_7 := $
\begin{minipage}{1.0in}\begin{center}
  \includegraphics{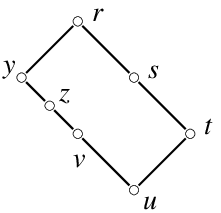}
\end{center}\end{minipage}
\hspace{0.2in}
$\Sc'_8 := $
\begin{minipage}{1.0in}\begin{center}
  \includegraphics{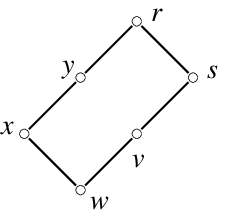}
\end{center}\end{minipage}
\hspace{0.1in} and \hspace{0.1in}
$\Sc'_9 := $
\begin{minipage}{1.0in}\begin{center}
  \includegraphics{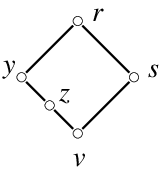}
\end{center}\end{minipage}
\end{center}
be cyclic subposets of $\Sc$.
Then, (vii) $\Sc'_7$ is neither an induced subposet nor
a subposet via induced cover subgraph,
(viii) $\Sc'_8$ is a subposet via induced cover subgraph but not an
induced subposet, and
(ix) $\Sc'_9$ is both an induced subposet and
a subposet via induced cover subgraph.
We note that if a cyclic subposet is an induced subposet,
then it must be a subposet via induced cover subgraph.
We will show this (Lemma~\THn{lem.subposet}) in Section~\SSn{cycle}.
\qed
\EndExample

\SubSSection[cycle]{Induced cyclic subposets}

\SStop
For developments later in this paper,
a study of cyclic subposets turns out to be crucial,
and for this we must also study path subposets.
Since the material here is irrelevant until Section~\SSn{class.b}, the reader may
wish to return to the present subsection after reading
Section~\SSn{me.summary}.

Let $\Sc$ be a poset.
A path or a cycle $\Vc$ (with ground set $V$) in $\Sc$ is by definition a subposet of
the subposet $\Vc'$ via induced cover subgraph of $\Sc$ on $V$,
and $\Vc'$ is in turn a subposet of the induced subposet $\Vc''$ of $\Sc$ on $V$.
So if this $\Vc$ is equal to $\Vc''$, then $\Vc = \Vc'$.
Thus, we have established

\Lemma[lem.subposet]
Let $\Vc$ be a path or a cycle in $\Sc$ and have ground set $V$.
If $\Vc$ is the induced subposet of $\Sc$ on $V$,
then $\Vc$ is the subposet via induced cover subgraph of $\Sc$ on $V$.
\EndLemma

An upward (or downward) path is a chain and therefore it
is both (i) an induced subposet and (ii) a subposet via induced cover
subgraph.
As shown in Example~\THn{ex.induce}(vi),
a path in general may be neither of these;
however, we can always devise a path with such properties which substitutes.

\Lemma[cons.path]
Suppose that there exists a path from $x$ to $y$ in $\Sc$.
Then there is a path $\Vc$ from $x$ to $y$ in $\Sc$ which
is an induced subposet of $\Sc$.
\EndLemma

\Proof
Partially ordering the up-sets and (separately) the down-sets in $\Sc$ by set inclusion,
let $U_0$ be a minimal up-set in $\Sc$ containing the vertices of some
path from $x$ to $y$.
(By assumption, $S$ is an up-set containing such
a path, so $U_0$ exists.)
Let $V_0$ be a minimal down-set in $\Sc$ such that $U_0 \cap V_0$
contains the vertices of some path from $x$ to $y$.
(Again, $S$ is a down-set satisfying this condition, so $V_0$ exists.)
Let $\Wc$ be a path from $x$ to $y$.
We can label minimal and maximal elements of $\Wc$
and count the segments of $\Wc$ alternating upward and downward as follows:
As the path $\Wc$ is traversed from $x$ to $y$, the path traces out either an upward
or a downward path from $z_0 = x$ to $z_1$,
either a downward or an upward path from $z_1$ to $z_2$, etc.,
alternatingly, as illustrated in
\begin{center}
$\Wc = $
\begin{minipage}{2.8in}\begin{center}
  \includegraphics{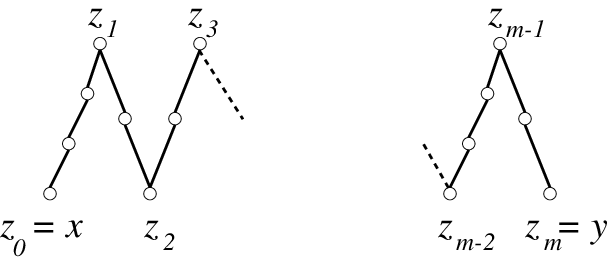}
\end{center}\end{minipage}
\end{center}
Then we can find a path $\Wc$ in $U_0\cap V_0$, with ground set $W$,
from $x$ to $y$ so that $m \ge 1$ is as small as possible,
that is, a path $\Wc$ from $x$ to $y$
satisfying $W \subseteq U_0 \cap V_0$ such that
the number $m$ of segments alternating upward and downward is smallest among such paths.

We first claim that any two minimal elements in $\Wc$ are
incomparable in $\Sc$ and that any two maximal elements in $\Wc$
are incomparable in $\Sc$.
To see this, suppose that $z_i > z_{j}$ for
two minimal elements $z_i$ and $z_{j}$ in $\Wc$.
Then there is some downward path $(u_0,\ldots,u_l)$ in $\Sc$ from $z_i = u_0$
to $u_l = z_{j}$.
Let $u_k < u_0$ be the first element among the $u_i$'s with $i \neq 0$
which belongs to $W$.
So $u_1,\ldots,u_{k-1} \not\in W$.
Since $u_0$, $u_k$ are in $U_0\cap V_0$, all the vertices of the
downward path $(u_0,\ldots,u_k)$ are also.
Replacing the part of $\Wc$ between $u_0$ and $u_k$ by the path 
$(u_0,\ldots,u_k)$ (possibly in reverse order), we get a path from $x$
to $y$ with a lower value of $m$, which is impossible.
Thus, any two minimal elements in $\Wc$ are incomparable in $\Sc$.
The same holds for any two maximal elements in $\Wc$.

To finish the proof, we claim that the path $\Wc$ is the induced subposet of $\Sc$ on $W$,
that is, that no pairs of elements of $W$ are comparable in $\Sc$,
beyond those specified by the poset $\Wc$.
To see this, suppose that there are elements $w$ and $w'$ of $W$ which are
comparable in $\Sc$ but not in $\Wc$.
By a replacement scheme similar to that employed in the preceding paragraph,
we can get a path $\Wc'$ (with ground
set $W'$) from $x$ to $y$ in $U_0 \cap V_0$ which bypasses some $z_i$.
If we define $V'_0 := \pideal{W'}$ to be the down-set in $\Sc$ generated
by $W'$ and (dually) $U'_0$
to be the up-set
$$
  U'_0 := \{\xi\in S: \eta\le\xi \text{ for some $\eta\in W'$}\},
$$
then $U'_0 \subseteq U_0$ and $V'_0 \subseteq V_0$, and
$W' \subseteq U'_0 \cap V'_0$.
If $z_i$ is minimal in $\Wc$, then
(by the preceding paragraph) $z_i$ is incomparable in $\Sc$
with any minimal element of $\Wc'$, and therefore $z_i \not\in U'_0$.
Similarly, if $z_i$ is maximal in $\Wc$, then $z_i \not\in V'_0$.
Thus, we have contradicted the minimality of $U_0$ or $V_0$
according as $z_i$ is minimal or maximal in $\Wc$,
which completes the proof.
\QED

The next three lemmas make it possible to construct an induced cyclic subposet
with specific properties when a poset $\Sc$ has a cycle.
The first of these lemmas is a simple corollary to
Lemma~\THn{cons.path} which ensures that
any non-acyclic poset has an induced cyclic subposet.

\Lemma[cons.cycle]
Suppose that a poset $\Sc$ has a cycle $(x_0,x_1,\ldots,x_{n-1},x_0)$.
Then $\Sc$ has an induced cyclic subposet
$(y_0,y_1,\ldots,y_{m-1},y_0)$ such that $x_0 = y_0$ and $x_{n-1} = y_{m-1}$.
\EndLemma

\Proof
Let $\Sc'$ be the poset on ground set $S$ obtained by deleting the
edge joining $x_{n-1}$ and $x_0$ from the Hasse diagram of $\Sc$.
Since $(x_0,x_1,\ldots,x_{n-1})$ is a path from $x_0$ to $x_{n-1}$
in $\Sc'$, by Lemma~\THn{cons.path} there is a
path $(y_0,y_1,\ldots,y_{m-1})$ from $y_0 = x_0$ to $y_{m-1} =
x_{n-1}$ which is an induced subposet of $\Sc'$.
Then, the cycle $(y_0,y_1,\ldots,y_{m-1},y_0)$ is as desired.
\QED

\Lemma[diamond]
If a poset $\Sc$ has a pair $(x,y)$ of elements such that
there exist at least two unequal upward paths from $x$ to $y$ in $\Sc$,
then $\Sc$ has an induced cyclic subposet which is a subdivided diamond.
\EndLemma

\Proof
Let $(u_0,\ldots,u_k)$ and $(v_0,\ldots,v_l)$ be two unequal
upward paths from $x$ to $y$ in $\Sc$.
Without loss of generality, the two paths differ in their first step,
i.e., $u_1 \neq v_1$.
Clearly $u_1$ and $v_1$ are incomparable.
Let $U := \{\xi\in S: u_1 < \xi \text{ and } v_1 < \xi\}$.
Then $U$ is nonempty because $y \in U$.
Let $y'$ be a minimal (in $\Sc$) element of $U$.
Let $W_1$ and $W_2$ be upward paths from $u_1$ to $y'$ and from $v_1$
to $y'$, respectively.
Then, for any $\xi\in W_1\setminus\{y'\}$ and $\eta\in W_2\setminus\{y'\}$,
$\xi$ and $\eta$ are incomparable; otherwise, the minimality of $y'$
is contradicted.
Thus, the ground set $\{x\}\cup W_1\cup W_2$ gives the desired induced subposet.
\QED

\Lemma[cycle.height]
If a poset $\Sc$ has a cycle with height at least $3$, then
$\Sc$ has an induced cyclic subposet $\Vc$ with height at
least $3$.
\EndLemma

\Proof
If $\Sc$ has a diamond as an induced subposet,
then by Lemma~\THn{diamond} we can find a subdivided diamond $\Vc$
which is an induced cyclic subposet of $\Sc$.
Clearly $\Vc$ has height at least $3$.

Suppose now that $\Sc$ has {\em no\/} diamond as an induced subposet.
Let $(x_0,x_1,\ldots,x_{n-1},x_0)$ be a cycle with height at least $3$
such that $x_{1} < x_0 < x_{n-1}$.
Let $S' := S \setminus\{x_0\}$ and let $\Sc'$ with ground set $S'$
be the subposet of $\Sc$ via induced cover graph.
So $(x_1,\ldots,x_{n-1})$ is a path in $\Sc'$.
By Lemma~\THn{cons.path}, there is a path
$\Uc = (u_0,u_1,\ldots,u_{k-1})$ in $\Sc'$
from $u_0 = x_1$ to $u_{k-1} = x_{n-1}$ which is an induced subposet
of $\Sc'$.
If $\Uc$ is an upward path from $x_1$ to $x_{n-1}$,
then $\Sc$ has two distinct upward paths from $x_1$ to $x_{n-1}$
which, by Lemma~\THn{diamond}, contradicts our assumption.
Thus, $\Uc$ is not an upward path and in particular
$x_1$ and $x_{n-1}$ are incomparable in $\Uc$ but comparable in
$\Sc$.
This implies that $\Uc$ is not an induced subposet of $\Sc$.
Let $\Uc' = (u_i,u_{i+1},\ldots,u_{i'})$ be a minimal segment of the path
$\Uc$ which is {\em not\/} an induced subposet of $\Sc$, that is,
a segment $\Uc'$ such that (i) $\Uc'$ is not an induced subposet of
$\Sc$, and (ii) any proper segment of $\Uc'$ is an induced subposet of
$\Sc$.
Then $u_i$ and $u_{i'}$ must be incomparable in $\Uc$ but comparable
in $\Sc$.
Without loss of generality, we may assume that $u_i < u_{i'}$ in $\Sc$.
Then there is a downward path $\Wc = (w_0,w_1,\ldots,w_{l-1})$ in $\Sc$
from $w_0 = u_{i'}$ to $w_{l-1} = u_{i}$.

Let $\Vc := (u_i,u_{i+1},\ldots,u_{i'},w_1,\ldots,w_{l-1})$ be a
cycle in $\Sc$.
Since $\Uc$ is an induced subposet of $\Sc'$, we must have $x_0 \in W$
and in particular $\Vc$ has height at least $3$.
We claim that $\Vc$ is an induced subposet of $\Sc$;
establishing the claim will complete the proof of the lemma.
To prove the claim, suppose that $u_j$ and $w_{j'}$ are comparable in $\Sc$
but incomparable in $\Vc$ for some $u_j \in U'\setminus\{u_i,u_{i'}\}$ and $w_{j'} \in W$.
Then $u_j$ is comparable with $u_i$ or $u_{i'}$ in $\Sc$ according as
$w_{j'} < u_j$ or $w_{j'} > u_j$.
If the pair (either $u_j,u_i$ or $u_j,u_{i'}$) are comparable in $\Vc$,
then there are two unequal upward paths (the one through $w_{j'}$ and
the other consisting of a segment of $\Vc$) with common ends in $\Sc$.
By Lemma~\THn{diamond}, this contradicts our diamond-free assumption.
If the pair is incomparable in $\Vc$,
then some proper segment of $\Uc'$ is not an induced subposet of $\Sc$,
and this contradicts the minimality of $\Uc'$.
\QED

Lemma~\THn{cons.cycle} implies that
if a poset $\Sc$ is non-acyclic then it has an induced cyclic subposet.
The next result gives various {\em sufficient\/} conditions for $\Sc$
to be non-acyclic.

\Proposition[non.acyclic]
If a poset $\Sc$ has an induced subposet $\Sc'$
poset-isomorphic to any of the following posets,
then $\Sc$ is non-acyclic:
\romanList
\item the diamond;
\item a subdivided crown with height at least $3$;
\item the $k$-crown for some $k \ge 3$;
\item the following ``double-bowtie'' poset:
\EndList
\Equation[s.dbow]
\begin{minipage}{1.2in}\begin{center}
  \includegraphics{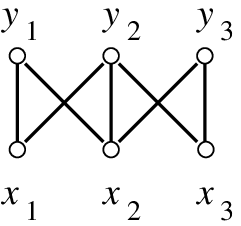}
\end{center}\end{minipage}
\end{equation}
\EndProposition

\Proof
(i) If $\Sc'$ is the diamond as in~\EQn{diamond},
then there are at least two unequal upward paths from $x$ to $w$.
By Lemma~\THn{diamond}, $\Sc$ has a cycle (namely, a subdivided diamond).

(ii) Suppose that $\Sc'$ is a subdivision of the $k$-crown
displayed and labeled in Figure~\ref{fig.posets}(d).
Since by assumption $\Sc'$ has height at least $3$,
without loss of generality
we may assume that there exists $z'\in S'$ such that $x_0 < z' < y_0$.
Then we can find an upward path $(x_0,\ldots,z,z',\ldots,y_0)$ in $\Sc$
from $x_0$ to $y_0$ with height at least $3$.
Since $z'$ is incomparable in $\Sc$ with each of $x_1,\ldots,x_{k-1}$ and
each of $y_1,\ldots,y_{k-1}$,
no upward path in $\Sc$ from any $x_i$ to any $y_j$ contains the directed edge $(z,z')$
unless $(i,j) = (0,0)$.
Let $(S,\Ec_{\Sc})$ be the cover graph of $\Sc$.
Then $x_0$ and $z'$ are connected in the graph
$(S,\Ec_{\Sc}\setminus\{\{z,z'\}\})$, which implies that
$\Sc$ is non-acyclic.

(iii) Suppose that $\Sc'$ is the $k$-crown for some $k \ge 3$, as
displayed and labeled in Figure~\ref{fig.posets}(d).
Since the set $B := \{\xi\in S: x_0\le\xi\le y_0
\text{ and } x_0\le\xi\le y_{k-1}\}$ is nonempty,
we can find a maximal element $x_0'$ in the set $B$.
Then $x_0'$ is incomparable with $x_1,\ldots,x_{k-1},y_1,\ldots,y_{k-2}$
and therefore the subposet of $\Sc$ induced by the ground set
$\{x_0',x_1,\ldots,x_{k-1},y_0,\ldots,y_{k-1}\}$ is again a $k$-crown.
Thus, we may without loss of generality
assume that the $k$-crown has no element $\xi$
satisfying $x_0 < \xi\le y_0$ and $x_0 < \xi\le y_{k-1}$.
Let $(x_0,z,\ldots,y_0)$ (with $z = y_0$ possible)
be an upward path from $x_0$ to $y_0$.
By our assumption, an upward path from $x_0$ to $y_{k-1}$ does not
contain the directed edge $(x_0,z)$.
Furthermore, no upward path from any $x_i$ to any $y_j$
[except when $(i,j) = (0,0)$] contains
the directed edge $(x_0,z)$ either;
otherwise, $\Sc'$ is not an induced subposet of $\Sc$.
Thus we see that $x_0$ and $z$ are connected in the graph
$(S,\Ec_{\Sc}\setminus\{\{x_0,z\}\})$,
which implies that $\Sc$ is non-acyclic.

(iv) Suppose that $\Sc'$ is the poset as in~\EQn{s.dbow}.
Since the set $B :=
\{\xi\in S: x_1\le\xi\le y_1 \text{ and } x_1\le\xi\le y_2\}$
is nonempty, we can find a maximal element $x_1'$ in the set $B$.
If $x_1'$ and $x_2$ are comparable, then we have $x_2 < x_1'$.
Noticing that $x_1'$ is incomparable with $x_3$ and $y_3$,
the cycle $(x_2,x_1',y_2,x_3,y_3,x_2)$ is an induced
subposet of $\Sc$ which is a subdivided $2$-crown with height $3$.
Then (ii) implies that $\Sc$ is non-acyclic.
So we may assume that $x_1'$ and $x_2$ are incomparable.
But then an upward path $(x_1',z,\ldots,y_1)$ in $\Sc$ from $x_1'$ to $y_1$
(with $z = y_1$ possible) does not share the directed edge $(x_1',z)$
with any upward path from $x_2$ to either $y_1$ or $y_2$.
Moreover, an upward path from $x_1'$ to $y_2$ does not contain the directed
edge $(x_1',z)$; otherwise, maximality of $x_1'$ is contradicted.
Thus, we see that $x_1'$ and $z$ are connected in the graph
$(S,\Ec_{\Sc}\setminus\{\{x_1',z\}\})$,
which implies that $\Sc$ is non-acyclic.
\QED

\SubSSection[strassen]{Extensibility and Strassen's theorem}

\SStop
Strassen's pioneering work~\cite{Strassen}
on the existence of probability measures with specified marginals
has been influential for the development of
the theory and applications of stochastic ordering (e.g.,~\cite{KKO,Ruschendorf,Shortt}).
We will treat briefly the general subject of
probability measures with specified marginals and
review some results essential for our later development.
Since we restrict attention to finite sets in the present paper,
some of the results presented here are greatly simplified
by our not needing to deal with topological and other technical matters.
(For an interesting review of the subject matter in a general
topological setting, see~\cite{Ruschendorf}.)

Let $A$ and $S$ be finite sets and
let $S^A$ be the collection of all functions
$\xx = (x_\alpha: \alpha\in A)$ from $A$ into $S$.
For $\xx\in S^A$ and $\alpha\in A$,
$\pi_\alpha(\xx)$ will denote the $\alpha$-coordinate of $\xx$.
Let $\alpha\in A$ be fixed.
Then $\pi_\alpha$,
the {\em $\alpha$-projection\/} from $S^A$ to the
$\alpha$-coordinate space $S$,
is a surjective map from $S^A$ to $S$.
Given a probability measure $Q$ on $S^A$, we define the probability
measure $Q\circ\pi_\alpha^{-1}$ on $S$ in the usual way via
$(Q\circ\pi_\alpha^{-1})(B) := Q(\pi_\alpha^{-1}(B))$
for any subset $B$ of $S$.

Consider the set of all signed measures on $S^A$ as a
normed vector space equipped with a suitable topology.
Strassen established the following theorem.

\Theorem[strassen]{\bf (Strassen~\cite{Strassen})}
Let $\Lambda$ be a nonempty convex closed subset of probability measures
on $S^A$ and let $(P_\alpha: \alpha\in A)$ be a system of probability
measures on $S$.
Then there exists a probability measure $Q \in \Lambda$ such that
\Equation[ext.marg]
  Q\circ\pi_\alpha^{-1} = P_\alpha
  \quad\text{\rm for every $\alpha\in A$}
\end{equation}
if and only if
\Equation[strassen]
  \sum_{\alpha\in A}\left(
    \sum_{\xi\in S} P_\alpha(\{\xi\}) f_\alpha(\xi) 
  \right)
  \le
  \sup\left\{
    \sum_{\xx\in S^A}
      Q(\{\xx\})
      \left(\sum_{\alpha\in A} f_\alpha\circ\pi_\alpha \right)(\xx)
    : Q \in \Lambda
  \right\}
\end{equation}
for any system $(f_\alpha: \alpha\in A)$ of real-valued functionals
on $S$.
\EndTheorem

Let $\Delta$ be a nonempty subset of $S^A$.
Then we say that a system $(P_\alpha: \alpha\in A)$ of probability
measures on $S$ is {\em extensible on\/} $\Delta$ if there
exists a probability measure $Q$ on $S^A$ satisfying~\EQn{ext.marg} and
\Equation[ext.mass]
  Q(\Delta) = 1.
\end{equation}
Let $\Lambda_\Delta$ be the set of all probability measures on
$S^A$ satisfying~\EQn{ext.mass}.
Clearly, $\Lambda_\Delta$ is nonempty, closed, and convex, so
Theorem~\THn{strassen} applies to it.
Observe that $\Lambda_\Delta$ is the convex hull of the set
$\{\delta_{\xx}: \xx\in\Delta\}$ where $\delta_{\xx}$ denotes the
point-mass probability at $\xx$.
Thus, the following theorem is a special case of Theorem~\THn{strassen}.

\Theorem[ext.thm]
Let $\Delta$ be a nonempty subset of $S^A$
and let $(P_\alpha: \alpha\in A)$ be a system of probability
measures on $S$.
Then $(P_\alpha: \alpha\in A)$ is extensible on $\Delta$ if and only if
\Equation[ext.thm]
  \sum_{\alpha\in A}\left(
    \sum_{\xi\in S} P_\alpha(\{\xi\}) f_\alpha(\xi)
  \right)
  \le
  \sup\left\{
    \sum_{\alpha\in A} f_\alpha\circ\pi_\alpha(\xx)
    : \xx \in \Delta
  \right\}
\end{equation}
for any system $(f_\alpha: \alpha\in A)$ of real-valued functionals
on $S$.
\EndTheorem

\Remark
Throughout this paper, we use the term ``system''
in place of ``family'' to refer to a
collection of probability measures, random variables, or real-valued
functionals.
When a partial ordering on the index set $A$ is introduced in the
later discussion, the usage becomes more appropriate.
We have co-opted the term ``extensibility'' for a system of probability measures
from a use by Vorob'ev~\cite{Vorobev} in a slightly different setting.
Vorob'ev's ``extensibility'' problem is now generally called the
{\em marginal problem\/} in the literature (e.g.~\cite{Ruschendorf,Shortt}).
\EndRemark

\SubSSection[me.problem]{The monotonicity equivalence problem}

\SStop
Since realizable monotonicity always implies stochastic monotonicity, 
the {\em monotonicity equivalence problem\/} for a given pair
$(\Ac,\Sc)$ of posets is to either verify or disprove Statement~\THn{state.me}.
\Statement[state.me]
For the given pair $(\Ac,\Sc)$ of posets,
every stochastically monotone system $(P_\alpha: \alpha\in A)$
of probability measures on $S$ is realizably monotone.
\EndStatement

We first formulate the monotonicity equivalence problem as a special case of 
the extensibility problem of Section~\SSn{strassen}.
Let $\Ac$ and $\Sc$ be finite posets.
We say that an element $\xx = (x_\alpha: \alpha\in A)$ of $S^A$ is {\em monotone\/}
if $x_\alpha \le x_\beta$ in $\Sc$ whenever $\alpha\le\beta$ in $\Ac$.
Define $\Delta$ to be the collection of all monotone elements of $S^A$.
Given a stochastically monotone system $(P_\alpha: \alpha\in A)$
of probability measures on $S$,
we say that a probability measure $Q$ on $S^A$
{\em realizes the monotonicity\/} if it satisfies~\EQn{ext.marg} and~\EQn{ext.mass}.
Observe that a system $(\XX_\alpha: \alpha\in A)$ of
$S$-valued random variables is merely an $S^A$-valued random
variable distributed as a probability measure $Q$ on $S^A$.
Clearly, the existence of $(\XX_\alpha: \alpha\in A)$
satisfying~\EQn{rm.mono}--\EQn{rm.marg} 
is equivalent to the existence of a probability measure $Q$ on $S^A$
satisfying~\EQn{ext.marg} and~\EQn{ext.mass}.
Thus, $(P_\alpha: \alpha\in A)$ is realizably monotone if and only if
it is extensible on $\Delta$.
This formulation establishes Theorem~\THn{ext.thm}
as a necessary and sufficient condition for realizable monotonicity.

We now present, without proof, some first simple results on
the monotonicity equivalence problem.
The upshot of these results is that we need only consider connected
posets in our investigation of monotonicity equivalence.

\Lemma[a.induce]
Suppose that $\Ac'$ is a (not necessarily induced) subposet of $\Ac$.
If $(P_\alpha: \alpha\in A)$ is realizably monotone, then
so is $(P_\alpha: \alpha\in A')$.
\EndLemma

\Lemma[s.induce]
Suppose that $\Sc'$ is an induced subposet of $\Sc$.
If monotonicity equivalence holds for $(\Ac,\Sc)$,
then it holds for $(\Ac,\Sc')$.
\EndLemma

\Proposition[a.connect]
Suppose that $\Ac$ is the disjoint union of nonempty posets $\Ac_1$ and $\Ac_2$.
Then monotonicity equivalence holds for $(\Ac,\Sc)$
if and only if it holds for both $(\Ac_1,\Sc)$ and $(\Ac_2,\Sc)$.
\EndProposition

\Proposition[s.connect]
Suppose that $\Sc$ is the disjoint union of nonempty posets $\Sc_1$
and $\Sc_2$.
Then monotonicity equivalence holds for $(\Ac,\Sc)$ if and only if
it holds for both $(\Ac,\Sc_1)$ and $(\Ac,\Sc_2)$.
\EndProposition

The next proposition is an immediate consequence of the
observation that the collection $\Delta$ of all monotone elements of
$S^A$ for $(\Ac,\Sc)$ is equal to the corresponding collection
for $(\Ac^*,\Sc^*)$.

\Proposition[dual]
Monotonicity equivalence holds for $(\Ac,\Sc)$ if and only if
it holds for $(\Ac^*,\Sc^*)$.
\EndProposition

\SSection[me.summary]{Subclasses of connected posets}

\SStop
As explained by Propositions~\THn{a.connect} and~\THn{s.connect},
we assume without further notice that $\Ac$ and $\Sc$ are connected
posets throughout the remainder of the paper.
We partition the collection of connected posets $\Sc$ into the
following four subclasses.
We say that
\abcList
\item
$\Sc$ is in {\em Class $\B$\/}, denoted $\Sc \in \B$, if
  $\Sc$ has either a cycle or an induced bowtie;
\item
$\Sc$ is in {\em Class $\Y$\/}, denoted $\Sc \in \Y$, if
  (i) $\Sc \not\in \B$, and (ii) $\Sc$ has an induced Y-poset;
\item
$\Sc$ is in {\em Class $\W$\/}, denoted $\Sc \in \W$, if
  (i) $\Sc \not\in \B\cup\Y$, and (ii) $\Sc$ has an induced W-poset;
\item
$\Sc$ is in {\em Class $\Z$\/}, denoted $\Sc \in \Z$, if
  $\Sc \not\in \B\cup\Y\cup\W$.
\EndList

Note that a poset $\Sc$ in Class~B may be acyclic.
For example, let
\begin{center}
$\Sc := $
\begin{minipage}{0.8in}\begin{center}
  \includegraphics{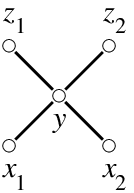}
\end{center}\end{minipage}
\end{center}
be an acyclic poset.
Then the subposet of $\Sc$ induced by $\{x_1,x_2,z_1,z_2\}$ is the
bowtie; thus, $\Sc \in \B$.
If the cover graph $(S,\Ec_{\Sc})$ of a given poset $\Sc$
has an element $x$ whose degree is at least $3$,
then $\Sc$ must have either a Y-poset or a W-poset as an
induced subposet and therefore $\Sc \in \B\cup\Y\cup\W$.
This implies that Class~Z consists precisely of those posets $\Sc$
whose cover graph $(S,\Ec_{\Sc})$ is a path
(and the nature of whose Hasse diagram is therefore ``zig-zag,'' which explains
our choice of ``Z'').

Given a poset $\Sc$, we call $\Ac$ a {\em poset of monotonicity
equivalence\/} or {\em of monotonicity inequivalence\/} (for $\Sc$)
according as Statement~\THn{state.me} is true or false for the pair $(\Ac,\Sc)$.
The question of monotonicity equivalence raised in Section~\SSn{intro}
can be recast as that of determining, for each $\Sc$, the class
$\Mc(\Sc)$ of all posets $\Ac$ of monotonicity equivalence for $\Sc$.
For a poset $\Sc$ in Class~B,~Y, or~Z, we can characterize the class
$\Mc(\Sc)$ precisely.
Furthermore, the class $\Mc(\Sc)$ is the {\em same\/} for
every $\Sc$ of the {\em same\/} class among Classes~B,~Y,~and~Z.
In the rest of this paper, we will show that
\begin{itemize}
\item
for every $\Sc\in B$,
$\Mc(\Sc)$ is the collection of all acyclic posets $\Ac$
(Theorem~\THn{class.b.main});
\item
for every $\Sc\in Y$,
$\Mc(\Sc)$ is the collection of posets $\Ac$ such that $\Ac$ is
enlargeable to an acyclic poset
(Theorem~\THn{class.y.main});
\item
for every $\Sc\in Z$,
$\Mc(\Sc)$ is the class of all posets $\Ac$
(Theorem~\THn{class.z.main}).
\end{itemize}

For a poset $\Sc$ of Class~W, we can exhibit a large subclass of
$\Mc(\Sc)$.
But the assertion that the class $\Mc(\Sc)$ is the same for
every $\Sc$ of Class~W is false.
Our investigation of Class~W is presented in the companion paper~\cite{invprob}.

\SSection[class.b]{The monotonicity equivalence problem on Class~B}

\SStop
In this section, we solve the monotonicity equivalence problem
when $\Sc$ is a poset of Class~B.
The main results of this section are summarized in the following two theorems.

\Theorem[a.acyclic]
If $\Ac$ is an acyclic poset, then monotonicity equivalence holds for
$(\Ac,\Sc)$ for any $\Sc$.
\EndTheorem

\Theorem[class.b.main]
Let $\Sc$ be a poset of Class~B.
Then monotonicity equivalence holds for $(\Ac,\Sc)$
if and only if $\Ac$ is an acyclic poset.
\EndTheorem

In Section~\SSn{a.acyclic} we briefly review a well-known
characterization of stochastic ordering and then prove Theorem~\THn{a.acyclic}.
Theorem~\THn{a.acyclic} establishes a sufficient condition for a poset $\Ac$
of monotonicity equivalence which is applicable to any poset $\Sc$.
But further generalization is not possible when $\Sc$ is a poset of Class~B.
In Section~\SSn{class.b.mi} we present various counterexamples where
monotonicity equivalence fails for a pair $(\Ac,\Sc)$ of non-acyclic posets.
In Section~\SSn{class.b.me}
we build on these counterexamples to
complete the proof of Theorem~\THn{class.b.main}.
Theorems~\THn{a.acyclic} and~\THn{class.b.main} can be immediately
combined to settle the monotonicity equivalence question
for Markov transition matrices, where $\Ac = \Sc$
(cf. the end of Section~\SSn{intro}):
\Theorem[me.markov]
If $\Ac = \Sc$, then
monotonicity equivalence holds for $(\Ac,\Sc)$
if and only if $\Sc$ is an acyclic poset.
\EndTheorem

\SubSSection[a.acyclic]{Stochastic ordering and acyclic index posets $\Ac$}

\SStop
To supplement the characterization of stochastic ordering
described in Section~\SSn{intro},
Kamae, Krengel, and O'Brien~\cite{KKO} introduced an equivalent
condition in terms of upward kernels.
Using their condition, they showed that,
for a sequence $(P_1,P_2,\ldots)$ of probability measures on a
common poset $\Sc$,
we have $P_i \stocle P_{i+1}$ for $i=1,2,\ldots$ if and only if
there exists a sequence $(\XX_1,\XX_2,\ldots)$
such that $\XX_i \le \XX_{i+1}$ and
$\Prob(\XX_i\in\cdot) = P_i(\cdot)$ for $i=1,2,\ldots$.
We will show (Theorem~\THn{a.acyclic}) that
this result can be generalized to a system $(P_\alpha: \alpha\in A)$
of probability measures whenever $\Ac$ is an acyclic poset.

Let $\Sc$ be a poset.
A function $k$ from $S \times 2^S$ to $[0,1]$ is
called a {\em stochastic kernel\/} on $S$ if
$k(x,\cdot)$ is a probability measure on $S$ for every $x\in S$.
A stochastic kernel $k$ on $\Sc$ is said to be {\em upward\/}
if
\Equation[up.kernel]
  k(x,\{\xi\in S: x \le \xi \text{ in $\Sc$}\}) = 1
  \quad\text{ for each $x \in S$. }
\end{equation}
We collect several characterizations of stochastic ordering in
the following proposition.

\Proposition[kko]{\bf (Kamae, Krengel, and O'Brien~\cite{KKO})}
Let $(P_1,P_2)$ be a pair of probability measures on $S$.
Then the following conditions are equivalent:
\abcList
\item $P_1 \stocle P_2$;
\item $P_1(V) \ge P_2(V)$ for every down-set $V$ in $\Sc$;
\item there exists a pair $(\XX_1, \XX_2)$ of $S$-valued random variables
  satisfying~\EQn{rv.le}--\EQn{rv.marg};
\item there exists an upward kernel $k$ such that
\EndList
\Equation[kko.kernel]
  P_2(\cdot) = \sum_{x\in S} P_1(\{x\}) k(x,\cdot).
\end{equation}
\EndProposition

Now consider the monotonicity equivalence problem.
The equivalence of (a) and (c) in Proposition~\THn{kko} can be extended
to equivalence between stochastic monotonicity and realizable monotonicity
in the case that $\Ac$ is acyclic.
The precise result has already been stated as Theorem~\THn{a.acyclic}.

\Proof[{\bf Proof of Theorem~\THn{a.acyclic}.}]
We prove the claim of Theorem~\THn{a.acyclic} by induction over the
cardinality of $A$.
The claim is vacuous when $|A| = 1$.
We now suppose that the claim is true for an acyclic poset $\Ac'$
when $|A'| = n-1$ for fixed $n \ge 2$,
and consider an acyclic poset $\Ac$ with cardinality $|A| = n$.
Let $a$ be a leaf in $\Ac$.
Without loss of generality, we assume that $a$ is maximal in $\Ac$;
thus, there is a unique element $b$ which is covered by $a$.
We consider the subposet $\Ac'$ of $\Ac$ induced by the ground set $A' = A\setminus\{a\}$.

Let $(P_\alpha: \alpha\in A)$ be a stochastically monotone
system of probability measures on $S$.
Then the subsystem $(P_\alpha: \alpha\in A')$ is stochastically monotone.
Since $\Ac'$ is an acyclic poset and $|A'| = n-1$,
by the induction hypothesis there exists a probability measure $Q'$ on $S^{A'}$
which realizes the monotonicity.
Since $P_{b} \stocle P_{a}$,
by Proposition~\THn{kko} there exists an upward kernel $k$
satisfying~\EQn{kko.kernel} for the pair 
$(P_{b},P_{a})$ of probability measures.
We can define a probability measure $Q$ on $S^{A}$ by
$$
  Q(\{\xx\}) := Q'(\{\pi_{A'}(\xx)\})\cdot
  k(\pi_{b}(\xx),\{\pi_{a}(\xx)\})
  \quad\text{ for $\xx\in S^{A}$,}
$$
where $\pi_{A'}$ denotes the projection from $S^{A}$ to $S^{A'}$ and
$\pi_\alpha$ denotes the $\alpha$-projection from $S^{A}$ to $S$
for each $\alpha \in A$.
In words, this says simply that we couple together the
probability measures $P_\alpha$, $\alpha\in A'$ using $Q'$ and then extend the
multivariate coupling to $P_a$ using the upward kernel $k$ from $P_b$
to $P_a$.
Observe that $Q'$ couples the probability measures
$(P_\alpha$: $\alpha\in A')$ correctly,
and that $a > \alpha \in A'$ in $\Ac$ implies $b \ge \alpha$ in $\Ac'$.
So coupling $P_a$ to $P_b$ correctly automatically couples $P_a$ to
each $P_\alpha$ ($\alpha \in A'$) correctly.
Thus, $Q$ realizes the monotonicity of $(P_\alpha: \alpha\in A)$,
and therefore the claim holds for $\Ac$.
\QED

\SubSSection[class.b.mi]{Monotonicity inequivalence on Class~B}

\SStop
The objective of this subsection is to present several examples of monotonicity
inequivalence.
To establish such an example we must exhibit a pair $(\Ac,\Sc)$ of posets
and a specific system $(P_\alpha: \alpha\in A)$ of probability
measures on $S$ which is stochastically but not realizably monotone.
We have already done this when both $\Ac$ and $\Sc$ are diamonds:
see Example~\THn{ex.dd}.
Our simple examples, including Example~\THn{ex.dd},
will serve as building blocks for more complex counterexamples
that establish quite general negative results.

\Example[ex.bd]
Let
\begin{center}
$\Ac := $
\begin{minipage}{1.0in}\begin{center}
  \includegraphics{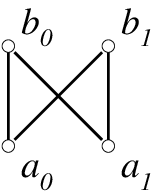}
\end{center}\end{minipage}
\end{center}
be the bowtie and let $\Sc$ be the diamond as in~\EQn{diamond}.
Define a system $(P_\alpha: \alpha\in A)$ of probability measures on
$S$ by
\Equation[ex.bd]
  P_\alpha := \begin{cases}
    \unif\{x,w\} & \quad\text{ if $\alpha = a_0$; } \\
    \unif\{y,z\} & \quad\text{ if $\alpha = a_1$; } \\
    \unif\{y,w\} & \quad\text{ if $\alpha = b_0$; } \\
    \unif\{z,w\} & \quad\text{ if $\alpha = b_1$. }
  \end{cases}
\end{equation}
The system is clearly stochastically monotone.

To see that it is not realizably monotone,
suppose that there exists a system $(\XX_\alpha: \alpha\in A)$ of
$S$-valued random variables which realizes the monotonicity.
Considering the event $\{\XX_{b_0} = y\}$, in order to maintain monotonicity
we must have
$$
  \Prob(\XX_{b_0} = y)
= \Prob(\XX_{a_0} = x,\:\XX_{a_1} = y,\:\XX_{b_0} = y,\:\XX_{b_1} = w)
= \frac{1}{2}.
$$
Similarly, we must have
$$
  \Prob(\XX_{a_0} = w)
= \Prob(\XX_{a_0} = w,\:\XX_{b_0} = w,\:\XX_{b_1} = w)
= \frac{1}{2}.
$$
Since the above two events are disjoint, we must have $\Prob(\XX_{b_1}=w) = 1$,
which contradicts $\Prob(\XX_{b_1}=w) = \frac{1}{2}$.
\qed
\EndExample

\Example[ex.bc]
Let $\Ac$ be the bowtie and let $\Sc$ be a $k$-crown with $k \ge 2$.
The posets $\Ac$ and $\Sc$ are displayed and labeled in
Example~\THn{ex.bd} and Figure~\ref{fig.posets}(d), respectively.
Define a system $(P_\alpha: \alpha\in A)$ of probability measures on
$S$ by
\Equation[ex.bc]
  P_\alpha := \begin{cases}
    \frac{k-1}{k}\,\unif(S_1 \setminus \{x_1\})
    + \frac{1}{k}\,\unif(\{y_0,x_{1}\})
    & \quad\text{ if $\alpha = a_0$; } \\
    \frac{1}{k}\,\unif(\{x_0\})
    + \frac{k-1}{k}\,\unif(S \setminus \{x_0,y_0\})
    & \quad\text{ if $\alpha = a_1$; } \\
    \frac{1}{k}\,\unif(\{y_{k-1}\})
    + \frac{k-1}{k}\,\unif(S \setminus \{x_0,y_{k-1}\})
    & \quad\text{ if $\alpha = b_0$; } \\
    \frac{k-1}{k}\,\unif(\{x_0\} \cup (S_2 \setminus \{y_0,y_{k-1}\}))
    + \frac{1}{k}\,\unif(\{y_0,y_{k-1}\})
    & \quad\text{ if $\alpha = b_1$, }
  \end{cases}
\end{equation}
where
$$
S_1 := \{x_0,x_1,\ldots,x_{k-1}\} \:\text{ and }\:
S_2 := \{y_0,y_1,\ldots,y_{k-1}\}.
$$
Then $(P_\alpha: \alpha\in A)$ is stochastically monotone.

Now let $\Delta$ be the collection of all monotone elements of $S^A$
and let
$$
  U_\alpha := \begin{cases}
    \{y_0\}
    & \quad\text{ if $\alpha = a_0$; } \\
    S_2 \setminus\{y_0\}
    & \quad\text{ if $\alpha = a_1$; } \\
    S_1 \setminus\{x_0\}
    & \quad\text{ if $\alpha = b_0$; } \\
    \{x_0\}
    & \quad\text{ if $\alpha = b_1$. }
  \end{cases}
$$
This builds a system $(I_{U_\alpha}: \alpha\in A)$ of
real-valued functions on $S$,
where $I_{U_\alpha}$ denotes the indicator
function of a subset $U_\alpha$ of $S$.
It is not hard to verify that
$$
  \sum_{\alpha\in A} \left(I_{U_\alpha}\circ\pi_\alpha\right)(\xx) \le 1
  \quad\text{ for any $\xx\in\Delta$. }
$$
Since
$$
  \sum_{\alpha\in A} P_\alpha(U_\alpha) = 1 + \frac{1}{2k},
$$
by Theorem~\THn{ext.thm} we have shown that
$(P_\alpha: \alpha\in A)$ is not realizably monotone.
\qed
\EndExample

\Remark[P0.P1]
The specific systems $(P_\alpha: \alpha\in A)$
of probability measures on $S$ presented in Examples~\THn{ex.bd} and~\THn{ex.bc}
will be used in our later discussions.
We further define probability measures $P_{\hat{0}}$ and $P_{\hat{1}}$
on $S$ for each example.
\romanList
\item
In Example~\THn{ex.bd}, let $P_{\hat{0}} := \delta_x$ and $P_{\hat{1}}
:= \delta_w$.
Clearly, $P_{\hat{0}} \stocle P \stocle P_{\hat{1}}$
for any probability measure $P$ on $S$.
\item
In Example~\THn{ex.bc},
let $P_{\hat{0}} := \unif(S_1)$ and $P_{\hat{1}} := \unif(S_2)$.
Then we have $P_{\hat{0}} \stocle P_\alpha \stocle P_{\hat{1}}$
for every $P_\alpha$ defined at~\EQn{ex.bc}.
\EndList
\EndRemark

\Example[ex.dc]
Let
\Equation[a.diamond]
\Ac :=
\begin{minipage}{0.8in}\begin{center}
  \includegraphics{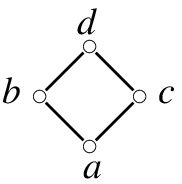}
\end{center}\end{minipage}
\end{equation}
be the diamond and
let $\Sc$ be a $k$-crown for $k \ge 2$
as in Figure~\ref{fig.posets}(d).
Define a system $(P_\alpha: \alpha\in A)$ of probability measures on
$S$ by
$$
  P_\alpha := \begin{cases}
    \unif(S_1)
    & \quad\text{ if $\alpha = a$; } \\
    \unif(\{y_0\}\cup(S_1 \setminus \{x_0\}))
    & \quad\text{ if $\alpha = b$; } \\
    \unif(\{y_{k-1}\}\cup(S_1 \setminus \{x_0\}))
    & \quad\text{ if $\alpha = c$; } \\
    \unif(S_2)
    & \quad\text{ if $\alpha = d$, }
  \end{cases}
$$
where $S_1$ and $S_2$ are defined as in Example~\THn{ex.bc}.
Then $(P_\alpha: \alpha\in A)$ is stochastically monotone.

Now let $\Delta$ be the collection of all monotone elements of $S^A$
and let
$$
  U_\alpha := \begin{cases}
    \{x_0\}
    & \quad\text{ if $\alpha = a$; } \\
    \{y_0\}\cup(S_1 \setminus \{x_0\})
    & \quad\text{ if $\alpha = b$; } \\
    \{y_{k-1}\}\cup(S_1 \setminus \{x_0\})
    & \quad\text{ if $\alpha = c$; } \\
    \emptyset
    & \quad\text{ if $\alpha = d$. }
  \end{cases}
$$
Then we have
$$
  \sum_{\alpha\in A} \left(I_{U_\alpha}\circ\pi_\alpha\right)(\xx) \le 2
  \quad\text{ for any $\xx\in\Delta$. }
$$
To see this, suppose that the sum is $3$ for some
monotone element $\xx$.
Then we must have
$(\pi_a(\xx),\pi_b(\xx),\pi_c(\xx)) = (x_0,y_0,y_{k-1})$, which is
impossible.
Since
$$
  \sum_{\alpha\in A} P_\alpha(U_\alpha) = 2 + \frac{1}{k},
$$
we deduce from Theorem~\THn{ext.thm} that
$(P_\alpha: \alpha\in A)$ is not realizably monotone.
\qed
\EndExample

Examples~\THn{ex.dd} and~\THn{ex.bd} both employ a certain
probabilistic argument which assumes the existence of certain random
variables and leads to a contradiction.
Here we introduce a lemma which is useful in conjunction with such
probabilistic arguments
when we extend monotonicity equivalence beyond our previously
considered counterexamples.

\Lemma[eq.dist.lem]
Let $\Sc$ be a poset and let $(\XX_1,\XX_2)$ be a pair of
$S$-valued random variables.
If $\Prob(\XX_1\in\cdot) = \Prob(\XX_2\in\cdot)$ and
$\Prob(\XX_1 \le \XX_2) = 1$, then $\Prob(\XX_1 = \XX_2) = 1$.
\EndLemma

\Proof
Notice that for any $\xi \in S$,
\begin{equation*}
  \Prob(\XX_1 \le \xi)
= \Prob(\XX_1 \le \xi,\: \XX_1 \le \XX_2)
\ge \Prob(\XX_2 \le \xi) +
    \Prob(\XX_1 = \xi < \XX_2).
\end{equation*}
Since $\Prob(\XX_1 \le \xi) = \Prob(\XX_2 \le \xi)$,
we deduce $\Prob(\XX_1 = \xi < \XX_2) = 0$.
Thus we obtain
\begin{equation*}
  \Prob(\XX_1 < \XX_2)
= \sum_{\xi\in S}\, \Prob(\XX_1 = \xi < \XX_2) = 0,
\end{equation*}
which completes the proof.
\QED

Now let
\Equation[a.crown]
\Ac_k :=
\begin{minipage}{2.4in}\begin{center}
  \includegraphics{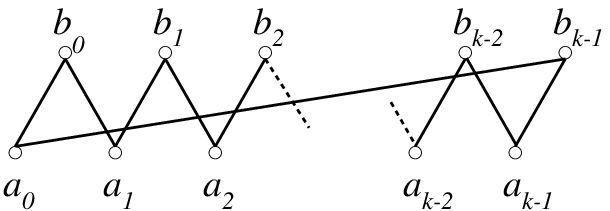}
\end{center}\end{minipage}
\end{equation}
be a $k$-crown.
If we have a known case of monotonicity inequivalence
for a pair $(\Ac_k,\Sc)$ of posets,
then we can apply Lemma~\THn{eq.dist.lem} to extend
monotonicity inequivalence to $(\Ac_{k'},\Sc)$ whenever $k' \ge k$.

\Proposition[a.crown.ext]
Let $\Ac_k$ be a $k$-crown as in~\EQn{a.crown}.
Given a pair $(\Ac_{k},\Sc)$ of posets,
suppose that there exists a stochastically monotone system $(P_\alpha: \alpha\in A_k)$ of
probability measures on $S$ which is not realizably monotone.
Then if $k' \ge k$, we can define
$$
  P_\alpha := \begin{cases}
    P_\alpha
    & \quad\text{\rm if $\alpha \in A_k$; } \\
    P_{b_{k-1}}
    & \quad\text{\rm if $\alpha \in A_{k'} \setminus A_k$, }
  \end{cases}
$$
to enlarge $(P_\alpha: \alpha\in A_k)$ to a stochastically monotone
system $(P_\alpha: \alpha\in A_{k'})$ which is not realizably monotone
for the pair $(\Ac_{k'},\Sc)$.
\EndProposition

\Proof
Since $P_{b_{k-1}} = P_{a_{k}} = \cdots = P_{a_{k'-1}} = P_{b_{k'-1}}$,
we see that the system $(P_\alpha: \alpha\in A_{k'})$ is stochastically monotone.
To see that it is not realizably monotone,
suppose that there exists a system $(\XX_\alpha: \alpha\in A_{k'})$ of
$S$-valued random variables which realizes the monotonicity
of $(P_\alpha: \alpha\in A_{k'})$.
By applying Lemma~\THn{eq.dist.lem} repeatedly,
we (almost surely) have
$\XX_{b_{k-1}} = \XX_{a_{k}} = \cdots = \XX_{a_{k'-1}} = \XX_{b_{k'-1}}$.
But then (after perhaps taking care of null sets)
$(\XX_\alpha: \alpha\in A_{k})$ realizes the monotonicity
of $(P_\alpha: \alpha\in A_{k})$, which is a contradiction.
\QED

As an immediate corollary to Proposition~\THn{a.crown.ext},
we can extend Examples~\THn{ex.bd} and~\THn{ex.bc}
to allow $\Ac$ to be the $k$-crown for arbitrary $k \ge 2$.
In summary, from the counterexamples in Example~\THn{ex.dd} and
Examples~\THn{ex.bd}--\THn{ex.dc} we have derived

\Proposition[mi.ex.thm]
Let $\Ac$ and $\Sc$ each be either a diamond or a crown.
Then monotonicity equivalence fails for $(\Ac,\Sc)$.
\EndProposition

\SubSSection[class.b.me]{The proof of Theorem~\THn{class.b.main}}

\SStop
Let $\Sc$ be a poset of Class~B.
Then we can find either (i) a $2$-crown as an induced subposet of $\Sc$
or (ii) a cycle as a (not necessarily induced) subposet of $\Sc$.
If $\Sc$ has a cycle, then,
by Lemma~\THn{cons.cycle}, $\Sc$ has an induced cyclic subposet $\Vc$.
It is possible to label the cycle $\Vc$
and to fix a starting point and orientation of the cycle so that,
as the cycle is traversed, it traces out an upward path from $z_0$
to $z_1$, then a downward path from $z_1$ to $z_2$,
then an upward path from $z_2$ to $z_3$,
etc., finishing with a downward path from $z_{2k-1}$ to $z_{0}$,
as illustrated in
\Equation[s.cycle]
\Vc =
\begin{minipage}{2.8in}\begin{center}
  \includegraphics{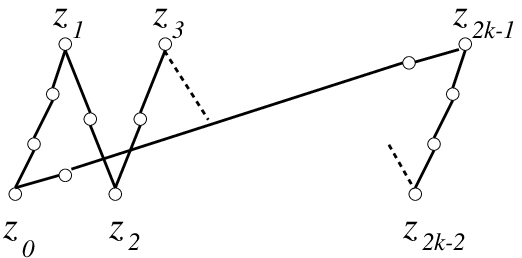}
\end{center}\end{minipage}
\end{equation}
If $k=1$, then $\Vc$ is a subdivided diamond;
otherwise, $\Vc$ is a subdivided $k$-crown ($k \ge 2$).
This observation gives a different characterization of Class~B.

\Lemma[class.b]
A poset $\Sc$ is in Class~B
if and only if $\Sc$ has either the diamond or a crown as an induced subposet.
\EndLemma

\Proof
We have already seen that a poset $\Sc$ of Class~B
has either the diamond or a crown as an induced subposet.
If $\Sc$ has an induced $2$-crown, then $\Sc$ is in Class~B by definition.
If $\Sc$ has an induced subposet which is
either the diamond or a $k$-crown for some $k \ge 3$,
then, by Proposition~\THn{non.acyclic}, $\Sc$ is non-acyclic and therefore
$\Sc$ is in Class~B, again by definition.
\QED

Now we turn to the proof of Theorem~\THn{class.b.main}.
If $\Ac$ is an acyclic poset, then,
by Theorem~\THn{a.acyclic}, monotonicity equivalence holds for
$(\Ac,\Sc)$.
Thus, the remaining task is to show that
if $\Ac$ is a non-acyclic poset, then monotonicity equivalence fails
for $(\Ac,\Sc)$.
By Lemmas~\THn{s.induce} and~\THn{class.b},
it suffices to show that monotonicity equivalence fails for $(\Ac,\Sc')$
whenever $\Sc'$ is either the diamond or a crown.
We complete the

\Proof[{\bf Proof of Theorem~\THn{class.b.main}.}]
Let $\Ac$ be a non-acyclic poset
and $\Sc'$ be either the diamond or an $m$-crown for some $m \ge 2$.
We will construct a stochastically monotone system $(P_\alpha: \alpha\in A)$
of probability measures on $S'$ which is not realizably monotone,
by dividing the construction into two cases.

\vspace{1.0ex}{\par\noindent\em Case~I.\/}
Suppose that $\Ac$ has a diamond $\Ac'$ as an induced subposet.
Let $\Ac'$ be labeled as in Example~\THn{ex.dc}.
By Examples~\THn{ex.dd} and~\THn{ex.dc},
there exists a stochastically monotone system $(P_\alpha: \alpha\in A')$
of probability measures on $S'$ which is not realizably monotone.
It then suffices by Lemma~\THn{a.induce} to show that
the system can be enlarged to a stochastically monotone system
$(P_\alpha: \alpha\in A)$ of probability measures on $S'$.

For this, define a partition $A_a,A_b,A_c$, and $A_d$ of $A$ by
$$
  A_\alpha := \begin{cases}
    A \setminus \{\alpha\in A: \text{ $b\le\alpha$ or $c\le\alpha$}\}
    & \quad\text{ if $\alpha = a$; } \\
    \{b\}
    & \quad\text{ if $\alpha = b$; } \\
    \{c\}
    & \quad\text{ if $\alpha = c$; } \\
    \{\alpha\in A: \text{ $b < \alpha$ or $c < \alpha$}\}
    & \quad\text{ if $\alpha = d$. }
  \end{cases}
$$
Then we can extend $(P_\alpha: \alpha\in A')$
to $(P_\alpha: \alpha\in A)$ by putting
$$
P_\alpha := P_\beta,\; \alpha\in A_\beta
$$
for each $\beta\in A'$.
It is routine to check that this extension maintains stochastic monotonicity,
that is, that if $\alpha_1 < \alpha_2$, then $P_{\alpha_1} \stocle P_{\alpha_2}$.
This is true if $\alpha_1,\alpha_2 \in A_\beta$ for some $\beta$
and also if $\alpha_1\in A_a$ or $\alpha_2\in A_d$.
If $\alpha_1 \in \{b,c\}$ and $\alpha_1 < \alpha_2$,
then $\alpha_2 \in A_d$.
If $\alpha_1 \in A_d$ and $\alpha_1 < \alpha_2$,
then $\alpha_2 \in A_d$.
So stochastic monotonicity is clear.

\vspace{1.0ex}{\par\noindent\em Case~II.\/}
Suppose that $\Ac$ has {\em no\/} diamond as an induced subposet.
By Lemma~\THn{cons.cycle}, $\Ac$ has an induced cyclic subposet $\Ac'$.
In the same way as what we did in~\EQn{s.cycle}, we can label the
cycle $\Ac'$ as illustrated in
\begin{center}
$\Ac' = $
\begin{minipage}{2.8in}\begin{center}
  \includegraphics{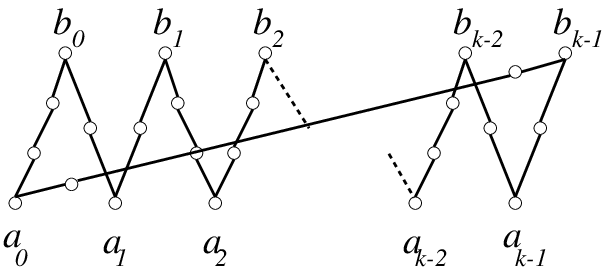}
\end{center}\end{minipage}
\end{center}
By our Case~II assumption, $\Ac'$ must be a subdivided $k$-crown for
some $k \ge 2$.
Let $\Ac''$ be the $k$-crown
$(a_0,b_0,a_1,b_1,\ldots,a_{k-1},b_{k-1},a_0)$.
By Example~\THn{ex.bd} and~\THn{ex.bc} (and then further using
Proposition~\THn{a.crown.ext}, if necessary),
there exists a stochastically monotone system $(P''_\alpha: \alpha\in A'')$ of
probability measures on $S$ which is not realizably monotone.

Let $P''_{\hat{0}}$ and $P''_{\hat{1}}$ be defined as in Remark~\THn{P0.P1}
so that $P''_{\hat{0}} \stocle P''_\alpha \stocle P''_{\hat{1}}$ for all
$\alpha\in A''$.
Consider the partition $\{A'_\beta$: $\beta\in A''\}$ of $A'$, where
$$
  A'_{a_i} := \{\alpha\in A':
    \text{ $a_i \le \alpha < b_{i-1}$ or $a_i \le \alpha < b_{i}$}\}
$$
and $A'_{b_i} := \{b_i\}$ for $i = 0,\ldots,k-1$.
By letting
$$
A_{\hat{1}} := \{\alpha\in A\setminus A': \alpha > \beta
\text{ for some $\beta\in A'$}\}
$$
and $A_{\hat{0}} := A \setminus (A'\cup A_{\hat{1}})$,
we can define a system $(P_\alpha: \alpha\in A)$ of
probability measures on $S$ by
$$
  P_\alpha := \begin{cases}
    P''_\beta
    & \quad\text{ if $\alpha \in A'_\beta$ for some $\beta\in A''$; } \\
    P''_{\hat{1}}
    & \quad\text{ if $\alpha \in A_{\hat{1}}$; } \\
    P''_{\hat{0}}
    & \quad\text{ if $\alpha \in A_{\hat{0}}$; }
  \end{cases}
$$
this system extends $(P''_\alpha: \alpha\in A'')$.

We claim that $(P_\alpha: \alpha\in A)$ is stochastically monotone.
Let $\alpha_1 < \alpha_2$.
If $\alpha_1\in A_{\hat{0}}$ or $\alpha_2\in A_{\hat{1}}$,
then $P_{\alpha_1} \stocle P_{\alpha_2}$.
This is also trivial if $\alpha_1,\alpha_2 \in A'$.
If $\alpha_1 \in A'$, then $\alpha_2 \in A'\cup A_{\hat{1}}$,
so $P_{\alpha_1} \stocle P_{\alpha_2}$.
If $\alpha_1 \in A_{\hat{1}}$, then $\alpha_2 \in A'\cup A_{\hat{1}}$.
We need only show that it is impossible to have both $\alpha_1 \in A_{\hat{1}}$
and $\alpha_2 \in A'$.
Indeed, then $\alpha_1 \not\in A'$, but for some $\beta\in A'$
we have $\beta < \alpha_1 < \alpha_2$.
But then there are two distinct upward paths from $\beta$ to $\alpha_2$
in $\Ac$, namely the one using edges in the cover graph $(A',\Ec_{\Ac'})$
and one containing $\alpha_1\not\in A'$.
This violates Lemma~\THn{diamond},
since we are assuming that $\Ac$ has no induced diamond.
Thus, we have established the claim and,
by Lemma~\THn{a.induce}, $(P_\alpha: \alpha\in A)$ cannot be
realizably monotone.
\QED

\SSection[class.y]{The monotonicity equivalence problem on Class~Y}

\SStop
In Section~\SSn{class.y} we investigate the
monotonicity equivalence problem when $\Sc \in \Y$.
The goal of this section is to prove the following theorem.
Two reformulations of the necessary and sufficient condition here are
given in Proposition~\THn{cyc.enlarge}.

\Theorem[class.y.main]
Let $\Sc$ be a poset of Class~$\Y$.
Then monotonicity equivalence holds for $(\Ac,\Sc)$ if and only if
there exists an acyclic poset $\tilde{\Ac}$ which has $\Ac$ as
an induced subposet.
\EndTheorem

Thus, some posets $\Ac$ of monotonicity equivalence may be non-acyclic.
As an instructive example, let
\Equation[bipart]
\Ac :=
\begin{minipage}{1.8in}\begin{center}
  \includegraphics{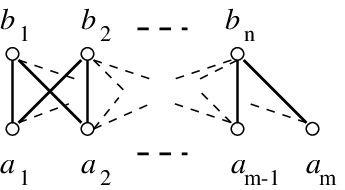}
\end{center}\end{minipage}
\end{equation}
be a poset where $a_i < b_j$ for all $i=1,\ldots,m$ and all
$j=1,\ldots,n$.
Then $\Ac$ is a poset of monotonicity equivalence for any $\Sc\not\in\B$.
To see this without resorting to Theorem~\THn{class.y.main},
grant the following proposition for now.

\Proposition[p.insert]
Suppose that $\Sc \not\in \B$.
Let $P_{a,1},\ldots,P_{a,m},P_{b,1},\ldots,P_{b,n}$ be probability measures
on $S$ satisfying
$$
  P_{a,i} \stocle P_{b,j}
  \quad\text{\rm for all $i=1,\ldots,m$ and all $j=1,\ldots,n$. }
$$
Then there exists a probability measure $P_0$ on $S$ such that
$$
  P_{a,i} \stocle P_0 \stocle P_{b,j}
  \quad\text{\rm for all $i=1,\ldots,m$ and all $j=1,\ldots,n$. }
$$
\EndProposition

Now suppose that $(P_\alpha: \alpha\in A)$ is a stochastically
monotone system of probability measures on $S$.
Proposition~\THn{p.insert} implies that there is a probability
measure $P_0$ on $S$ such that
$P_{a_i} \stocle P_0 \stocle P_{b_j}$
for all $i=1,\ldots,m$ and all $j=1,\ldots,n$.
Define an acyclic poset $\tilde{\Ac}$ on the set
$\tilde{A} = A\cup\{c\}$ by means of the Hasse diagram
\Equation[ext.bipart]
\tilde{\Ac} :=
\begin{minipage}{1.8in}\begin{center}
  \includegraphics{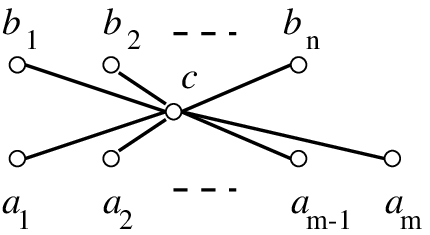}
\end{center}\end{minipage}
\end{equation}
Then the poset $\Ac$ is an induced subposet of $\tilde{\Ac}$.
By letting $P_{c} := P_0$, we can enlarge the system
$(P_\alpha: \alpha\in A)$ to the system $(P_\alpha: \alpha\in \tilde{A})$,
which remains stochastically monotone.
By applying Theorem~\THn{a.acyclic} and then Lemma~\THn{a.induce},
we see that both $(P_\alpha: \alpha\in \tilde{A})$ and $(P_\alpha: \alpha\in A)$
are realizably monotone.
Thus, we have shown that monotonicity equivalence holds for $(\Ac,\Sc)$.

In Section~\SSn{prob.on.acyc}
we attend to the proof of Proposition~\THn{p.insert}.
A large class of posets $\Ac$ of monotonicity inequivalence
is presented in Section~\SSn{class.y.mi},
leading to the proof of Theorem~\THn{class.y.main}
in Section~\SSn{class.y.me}.

\SubSSection[prob.on.acyc]{Probability measures on an acyclic poset}

\SStop
The goal of this subsection is to prove Proposition~\THn{p.insert}.
We begin this subsection by introducing a natural partial ordering
on a connected acyclic graph (i.e., a tree) when
one vertex is specified to become a top element,
that is, to be made larger than every other vertex.
Let $\Sc$ be a connected acyclic poset
and let $\tau$ be a fixed leaf of $\Sc$.
Declare $x \le_{\tau} y$ for $x,y \in S$
if and only if the (necessarily existent and unique)
path $(\tau,\ldots,x)$ from $\tau$ to $x$ contains
the path $(\tau,\ldots,y)$ from $\tau$ to $y$ as a segment.
This introduces another partial ordering $\le_{\tau}$ on the same ground set
$S$ (see~\cite{Bogart}).
We call this new poset $(S,\le_{\tau})$ a {\em rooted tree\/}
(rooted at $\tau$).
(Comparison of the poset $\Sc$ and a rooted tree is illustrated in
Example~\THn{ex.root.tree}.)
The element $\tau$ is clearly the maximum of the rooted tree $(S,\le_{\tau})$ and is
called the {\em root\/}.
If $x$ covers $y$ in $(S,\le_{\tau})$,
then $y$ is called a {\em successor\/} of $x$,
and $x$ is called the {\em predecessor\/} of $y$.

For each $x \in S$,
we define a {\em section of rooted tree\/} by
$$
  \csection{x} := \{\xi\in S: \xi \le_{\tau} x\},
$$
that is, the down set in $(S,\le_{\tau})$ generated by $x$
[cf. Section~\SSn{poset}(1)].
Every section $\csection{x}$ is either a down-set or an up-set in
$\Sc$, and which of these holds can be determined
from the cover relation of $\Sc$.
We state this as the following lemma.

\Lemma[lem.csec]
Let $\Sc$ be a connected acyclic poset.
For every $x \in S$,
$\csection{x}$ is either a down-set or an up-set in $\Sc$.
If $x \neq \tau$, then there is a unique predecessor $w$ of $x$,
and the edge $\{x,w\}$ belongs to the cover graph $(S,\Ec_{\Sc})$ of $\Sc$.
Moreover, $\csection{x}$ is {\rm (i)} a down-set or
{\rm (ii)} an up-set in $\Sc$ according as
{\rm (i)} $w$ covers $x$ or {\rm (ii)} $x$ covers $w$ in $\Sc$.
\EndLemma

\Proof
If $x=\tau$, then $\csection{\tau} = S$ is both a down-set and an
up-set in $\Sc$.
If $x \neq \tau$, then
$x <_{\tau} \tau$ and there is a unique predecessor of $x$;
otherwise, the uniqueness of the path is contradicted.
Let $w$ be the predecessor of $x$.
Clearly, $\{x,w\}$ belongs to the cover graph of $\Sc$.
Suppose that $w$ covers $x$ in $\Sc$.
We claim that $\csection{x}$ is a down-set in $\Sc$,
that is, that $\eta\in\csection{x}$ whenever 
$\eta \le \xi$ in $\Sc$ for some $\xi\in\csection{x}$.
[Since we have the same rooted tree $(S,\le_{\tau})$ for the dual $\Sc^*$,
in proving the claim we will also settle the case that $x$ covers $w$ in $\Sc$.]
To see this, look at the paths from the root $\tau$ to $\xi$ and $\eta$,
say, $(u_0,\ldots,u_{n-1})$ from $\tau = u_0$ to $u_{n-1} = \xi$
and $(v_0,\ldots,v_{m-1})$ from $\tau = v_0$ to $v_{m-1} = \eta$.
For some $k$, the two paths descend the same vertices until the $k$th vertex,
then split at the $(k+1)$st vertex.
The path from $\xi$ to $\eta$ is then
$(u_{n-1},\ldots,u_{k+1},u_k,v_{k+1},\ldots,v_{m-1})$,
which is downward in $\Sc$ by assumption.
That $\xi \le_{\tau} x$ implies that $(u_{i-1}, u_i) = (w,x)$ for some
$i$.
Furthermore, we have $i \le k$;
otherwise, the downward path from $\xi$ to $\eta$ contains the directed edge
$(x,w)$, which is impossible.
Thus, the path from $\tau$ to $\eta$ contains the vertex $x$,
which implies that $\eta \le_{\tau} x$.
\QED

\begin{figure}[h]
\vspace{0.2in}
\begin{center}
\begin{tabular}{ccc}
$\Sc :=$
\begin{minipage}{1.4in}\begin{center}
  \includegraphics{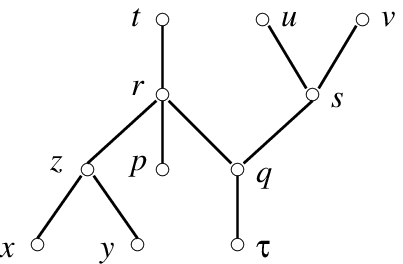}
\end{center}\end{minipage}
& \hspace{0.2in} &
$(S,\le_{\tau}) :=$
\begin{minipage}{1.6in}\begin{center}
  \includegraphics{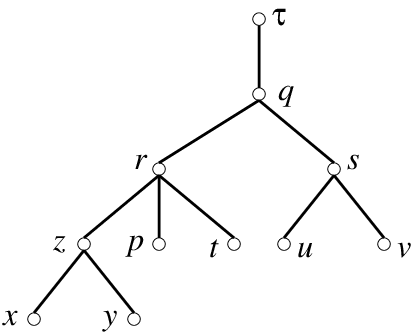}
\end{center}\end{minipage} \\[0.7in]
(a) a poset of Class~Y & & (b) a rooted tree with root $\tau$
\end{tabular}
\caption{The comparison of two posets}
\label{fig.root.tree}
\end{center}
\end{figure}

\Example[ex.root.tree]
Let $\Sc$ be the poset of Class~Y displayed in
Figure~\ref{fig.root.tree}(a).
By choosing the leaf $\tau$ of $\Sc$ as the root,
we obtain the rooted tree $(S,\le_{\tau})$
illustrated in Figure~\ref{fig.root.tree}(b).
For example, $r$ covers its predecessor $q$ in $\Sc$.
By Lemma~\THn{lem.csec},
the section $\csection{r} = \{r,p,t,x,y,z\}$ is an up-set in $\Sc$,
which we can confirm immediately from Figure~\ref{fig.root.tree}(a).
\qed
\EndExample

Now let $P$ be a probability measure on $S$.
We define the {\em distribution function\/} of $P$ by
\Equation[def.d.f]
  F(x) := P(\csection{x})
  \quad\text{ for each $x \in S$.}
\end{equation}
It satisfies
\Equation[pr.dist:1]
  F(\tau) = 1
\end{equation}
and
\Equation[pr.dist:2]
  \sum_{\xi\in\Cc(x)} F(\xi) \le F(x)
  \quad\text{ for every $x \in S$,}
\end{equation}
where $\Cc(x)$ denotes the set of all successors of $x$
(and the summation is defined to be zero if $\Cc(x) = \emptyset$).
Conversely, if a nonnegative function $F$ on $S$ satisfies
the properties~\EQn{pr.dist:1}--\EQn{pr.dist:2},
then it is the distribution function of the probability measure $P$
determined uniquely via
\Equation[pr.dist.ext]
  P(\{x\}) := F(x) - \sum_{\xi\in\Cc(x)} F(\xi)
  \quad\text{ for each $x \in S$.}
\end{equation}
Furthermore, stochastic ordering can be characterized
in terms of distribution functions, as stated in the following lemma.

\Lemma[csec.stoc]
Let $P_i$ be a probability measure on $S$
and let $F_i$ be the distribution function of $P_i$, for each $i=1,2$.
Then $P_1 \stocle P_2$ if and only if
for every $x\in S$ we have
\romanList
\item
$F_1(x) \le F_2(x)$
if $\csection{x}$ is an up-set in $\Sc$, and
\item
$F_1(x) \ge F_2(x)$
if $\csection{x}$ is a down-set in $\Sc$.
\EndList
\EndLemma

\Proof
By~\EQn{stocle.def} and its trivial consequence
Proposition~\THn{kko}(b),
$P_1 \stocle P_2$ clearly implies the conditions~(i)--(ii).
We proceed to the converse.
Since any up-set $U$ in $\Sc$ is the disjoint union of
the components $V_1,\ldots,V_m$ of the subgraph of
$(S,\Ec_{\Sc})$ induced by $U$ and $V_1,\ldots,V_m$ are all up-sets in
$\Sc$,
to prove $P_1 \stocle P_2$ it suffices to show~\EQn{stocle.def}
for every up-set $U$ which induces a connected subgraph of
$(S,\Ec_{\Sc})$.
If a set $U$ induces a connected subgraph of $(S,\Ec_{\Sc})$,
then we can find $x \in U$ and incomparable elements $y_1,\ldots,y_k$
of $\csection{x}$ in $(S,\le_{\tau})$ so that
$$
  U = \csection{x} \setminus \left(\bigcup_{i=1}^k \csection{y_i}\right),
$$
where (as usual) the union is empty if $k=0$.
Furthermore, suppose that $U$ is an up-set in $\Sc$.
If $x = \tau$, then $\csection{x} = S$ is trivially an up-set in
$\Sc$;
otherwise, $x$ covers its predecessor $w$ in $\Sc$
and, by Lemma~\THn{lem.csec}, $\csection{x}$ is an up-set in $\Sc$.
Similarly, $\csection{y_i}$ is a down-set in $\Sc$ for each $i=1,\ldots,k$.
Therefore, we have
$$
  P_1(U) = F_1(x) - \sum_{i=1}^k F_1(y_i)
  \le F_2(x) - \sum_{i=1}^k F_2(y_i) = P_2(U),
$$
which establishes the sufficiency of the conditions~(i)--(ii).
\QED

Because of Lemma~\THn{csec.stoc}, we write $F_1 \stocle F_2$ if a pair
$(F_1,F_2)$ of distribution functions on $S$ satisfies
Lemma~\THn{csec.stoc}(i)--(ii) for every $x \in S$.

We now turn to the proof of Proposition~\THn{p.insert}.
Let $\Sc\not\in\B$ and
let $F_{a,1},\ldots,F_{a,m}$, $F_{b,1},\ldots,F_{b,m}$ be the
distribution functions satisfying
$$
  F_{a,i} \stocle F_{b,j}
  \quad\text{ for all $i=1,\ldots,m$ and all $j=1,\ldots,n$. }
$$
Then define the function $\theta$ on $S$ by
\Equation[def.theta]
  \theta(x) := \begin{cases}
    \max\{ F_{a,i}(x): i=1,\ldots,m \}
    & \quad\text{ if $\csection{x}$ is an up-set in $\Sc$; } \\
    \max\{ F_{b,j}(x): j=1,\ldots,n \}
    & \quad\text{ if $\csection{x}$ is a down-set in $\Sc$. }
  \end{cases}
\end{equation}
for $x \in S$.
We first present the following lemma.

\Lemma[p.insert.lem]
Let $\Sc\not\in\B$.
Suppose that $x \in S$ and that $v_1,\ldots,v_l$
are mutually incomparable elements of $\csection{x}$ in $(S,\le_{\tau})$.
\abcList
\item If $\csection{x}$ is a down-set in $\Sc$ and $v_1,\ldots,v_l \le
x$ in $\Sc$, then
$$
\sum_{i=1}^l \theta(v_i) \le F_{a,j}(x)
\quad\text{\rm for all $j = 1,\ldots,m$.}
$$
\item
If $\csection{x}$ is an up-set in $\Sc$ and $v_1,\ldots,v_l \ge x$ in
$\Sc$, then
$$
\sum_{i=1}^l \theta(v_i) \le F_{b,j'}(x)
\quad\text{\rm for all $j' = 1,\ldots,n$.}
$$
\EndList
\EndLemma

\Proof
Suppose that the hypotheses in (a) hold.
If $v_1 = x$, then $l = 1$ and the inequality clearly holds.
Otherwise, $v_i \neq x$ for every $i=1,\ldots,l$.
Since the path $(v_i,u_i,\ldots,x)$ is upward and $u_i$ covers $v_i$ in $\Sc$,
by Lemma~\THn{lem.csec}(b) $\csection{v_i}$ is a down-set in $\Sc$.
Therefore we have
$$
  \sum_{i=1}^l \theta(v_i)
  \le \sum_{i=1}^l F_{a,j}(v_i) \le F_{a,j}(x)
$$
for all $j = 1,\ldots,m$, as desired.
The case (b) is reduced to (a) by considering the dual $\Sc^*$.
\QED

We now define a nonnegative function $F_0$ on $S$ inductively.
If $x$ is a minimal element in $(S,\le_{\tau})$, then
assign $F_0(x) := \theta(x)$.
If $x$ is a nonminimal element in $(S,\le_{\tau})$
and $F_0(\xi)$, $\xi\in\Cc(x)$, have all been assigned,
then set
$$
  F_0(x)
  := \max\left\{ \theta(x), \sum_{\xi\in\Cc(x)} F_0(\xi) \right\}.
$$
Clearly $F_0$ satisfies~\EQn{pr.dist:2}.
We complete the proof of Proposition~\THn{p.insert} by showing that
$F_0$ satisfies~\EQn{pr.dist:1} and
\Equation[p.insert.dist]
  F_{a,i} \stocle F_0 \stocle F_{b,j}
  \quad\text{ for all $i=1,\ldots,m$ and all $j=1,\ldots,n$. }
\end{equation}
Thus, $F_0$ is a distribution function with the desired property.

\Proof[{\bf Proof of Proposition~\THn{p.insert}.}]
We first claim that for every $x \in S$, there are incomparable elements
$v_1,\ldots,v_l$ of $\csection{x}$ in $(S,\le_{\tau})$
which satisfy both the hypotheses of one of
Lemma~\THn{p.insert.lem}(a),(b) and also
\Equation[p0.theta]
  F_0(x) = \sum_{i=1}^l \theta(v_i).
\end{equation}
We will show this by induction over the cardinality of $\csection{x}$.
If $|\csection{x}| = 1$, then $x$ is a minimal element in $(S,\le_{\tau})$
and indeed $F_0(x) = \theta(x)$.

Suppose that the claim holds for any $x \in S$ such that
$|\csection{x}| \le n-1$.
Let $ x\in S$ satisfy $|\csection{x}| = n \ge 2$.
If $x = \tau$, then recall that $\tau$ is a leaf in $\Sc$ so that $\Cc(\tau)$ is a
singleton, say $\{y\}$.
By the induction hypothesis, we can find incomparable elements
$v_1,\ldots,v_{l}$ of $\csection{y}$
in $(S,\le_{\tau})$ which satisfy both the hypotheses of
one of Lemma~\THn{p.insert.lem}(a),(b) and ~\EQn{p0.theta} for $y$.
If Lemma~\THn{p.insert.lem}(a) obtains, then
$$
  F_0(y) = \sum_{j=1}^{l} \theta(v_j)
  \le F_{a,1}(y) \le 1.
$$
A similar derivation concludes that $F_0(y) \le 1$
when Lemma~\THn{p.insert.lem}(b) obtains.
Therefore, we have $F_0(\tau) = \theta(\tau) = 1$,
which proves~\EQn{pr.dist:1}.
Furthermore, the claim holds for $x = \tau$.

If $x \neq \tau$, then $x$ has a predecessor $y_0$ and successors $y_1,\ldots,y_k$
for some $k \ge 1$ (recalling our assumption $|\csection{x}| = n \ge 2$).
Without loss of generality, we may assume that $y_0$ covers $x$ in $\Sc$.
(We can treat the case that $x$ covers $y_0$ in $\Sc$ in exactly the
same way by considering the dual $\Sc^*$.)
Then, by Lemma~\THn{lem.csec}, $\csection{x}$ is a down-set in $\Sc$.
Since $\Sc$ cannot have a bowtie as an induced subposet,
only the following three cases can occur:
(I) $x$ covers $y_1,\ldots,y_k$ in $\Sc$, or
(II) $y_1,\ldots,y_k$ cover $x$ in $\Sc$, or,
with $k \ge 2$, (III) $y_1,\ldots,y_{k-1}$ cover $x$, and $x$ covers $y_k$, in $\Sc$.
The induced subposet of $\Sc$ on $\{x,y_0,y_1,\ldots,y_k\}$ for each
of these three cases is illustrated in the following figure.
\begin{center}
\begin{tabular}{ccccc}
\begin{minipage}{1.0in}\begin{center}
  \includegraphics{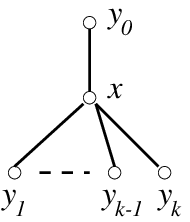}
\end{center}\end{minipage}
& \hspace{0.4in} &
\begin{minipage}{1.0in}\begin{center}
  \includegraphics{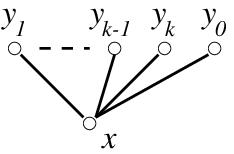}
\end{center}\end{minipage}
& \hspace{0.4in} &
\begin{minipage}{1.0in}\begin{center}
  \includegraphics{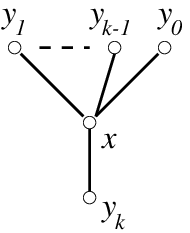}
\end{center}\end{minipage} \\[0.3in]
Case I & & Case II & & Case III
\end{tabular}
\end{center}

\vspace{1.0ex}{\par\noindent\em Case~I.\/}
If $F_0(x) = \theta(x)$, the claim is obvious.
Otherwise, we have $F_0(x) = \sum_{i=1}^k F_0(y_i)$.
For each $i=1,\ldots,k$, the section $\csection{y_i}$ is a down-set in $\Sc$
by Lemma~\THn{lem.csec}, and therefore
by the induction hypothesis we have
incomparable elements $v^{(i)}_1,\ldots,v^{(i)}_{l_i}$
of $\csection{y_i}$ in $(S,\le_{\tau})$
satisfying Lemma~\THn{p.insert.lem}(a) and~\EQn{p0.theta} for $y_i$.
Thus we have
$$
  F_0(x) = \sum_{i=1}^k F_0(y_i)
  = \sum_{i=1}^k \sum_{j=1}^{l_i} \theta(v^{(i)}_j),
$$
and $x \ge v^{(i)}_j$ for all $i,j$.
Since the $v^{(i)}_j$'s are incomparable in $(S,\le_{\tau})$,
the claim holds for $x$.

\vspace{1.0ex}{\par\noindent\em Case~II.\/}
For each $i=1,\ldots,k$, $\csection{y_i}$ is an up-set in $\Sc$
by Lemma~\THn{lem.csec}, and therefore
by the induction hypothesis we have
incomparable elements $v^{(i)}_1,\ldots,v^{(i)}_{l_i}$
of $\csection{y_i}$ in $(S,\le_{\tau})$
satisfying Lemma~\THn{p.insert.lem}(b) and~\EQn{p0.theta} for $y_i$.
By applying Lemma~\THn{p.insert.lem}(b) to~\EQn{p0.theta}, we have
$$
  \sum_{i=1}^k F_0(y_i)
  =   \sum_{i=1}^k \sum_{j=1}^{l_i} \theta(v^{(i)}_j)
  \le \sum_{i=1}^k F_{b,1}(y_i)
  \le F_{b,1}(x) \le \theta(x),
$$
which implies that $F_0(x) = \theta(x)$.
Thus the claim holds for $x$.

\vspace{1.0ex}{\par\noindent\em Case~III.\/}
By Lemma~\THn{lem.csec}, $\csection{y_k}$ is a down-set in $\Sc$
and therefore by the induction hypothesis
we can find incomparable elements
$v^{(k)}_1,\ldots,v^{(k)}_{l_k}$ of $\csection{y_k}$ in $(S,\le_{\tau})$
satisfying Lemma~\THn{p.insert.lem}(a) and~\EQn{p0.theta} for $y_k$.
Since $\Sc$ has no bowtie as an induced subposet
and $y_0,y_1 \ge y_k \ge v^{(k)}_1,\ldots,v^{(k)}_{l_k}$ in $\Sc$,
we have $l_k = 1$.
Thus, we can find some $j_0$ so that
$$
  F_0(y_k) = \theta(v^{(k)}_1) = F_{b,j_0}(v^{(k)}_1)
  \le F_{b,j_0}(y_k).
$$
For each $i=1,\ldots,k-1$,
the section $\csection{y_i}$ is an up-set in $\Sc$ by Lemma~\THn{lem.csec},
and therefore by the induction hypothesis we have incomparable elements
$v^{(i)}_1,\ldots,v^{(i)}_{l_i}$ of $\csection{y_i}$ in $(S,\le_{\tau})$
satisfying Lemma~\THn{p.insert.lem}(b) and~\EQn{p0.theta} for $y_i$.
By applying Lemma~\THn{p.insert.lem}(b) to~\EQn{p0.theta}, we obtain
\begin{align*}
  \sum_{i=1}^k F_0(y_i)
& \le F_{b,j_0}(y_k)
      + \sum_{i=1}^{k-1} \sum_{j=1}^{l_i} \theta(v^{(i)}_j) \\
& \le \sum_{i=1}^k F_{b,j_0}(y_i)
  \le F_{b,j_0}(x) \le \theta(x),
\end{align*}
which implies that $F_0(x) = \theta(x)$.
Thus, we have established the claim.

In order to show~\EQn{p.insert.dist}, it suffices to show that
if $\csection{x}$ is an up-set in $\Sc$ then we have
$$
  F_{a,i}(x) \le F_0(x) \le F_{b,j}(x)
  \quad\text{ for all $i=1,\ldots,m$ and all $j=1,\ldots,n$. }
$$
(Again, we can verify the case that $\csection{x}$ is a down-set in $\Sc$
by considering the dual $\Sc^*$.)
Suppose that $\csection{x}$ is an up-set in $\Sc$.
Then we can find incomparable elements
$v_1,\ldots,v_l$ of $\csection{x}$ in $(S,\le_{\tau})$
satisfying Lemma~\THn{p.insert.lem}(b) and~\EQn{p0.theta}.
By applying Lemma~\THn{p.insert.lem}(b) to~\EQn{p0.theta}, we have
$$
  F_{a,i}(x) \le \theta(x) \le F_0(x)
  = \sum_{i=1}^l \theta(v_i)
  \le F_{b,j}(x).
$$
This completes the proof.
\QED

\SubSSection[class.y.mi]{Monotonicity inequivalence on Class~Y}

\SStop
In this subsection, we present various examples,
each with a poset $\Sc$ from Class~Y, of posets $\Ac$ of monotonicity inequivalence.
The next example turns out to be a building block for all the other examples.

\Example[dia.y.ex]
Let $\Ac_0$ be the diamond given in~\EQn{a.diamond} and
let $\Sc_0$ be the Y-poset as in Figure~\ref{fig.posets}(b).
Define a system $(P_a,P_b,P_c,P_d)$ of probability measures on $S_0$ by
\Equation[dia.y]
P_\alpha := \begin{cases}
  \unif\{x,y\} & \quad\text{ if $\alpha = a$; } \\
  \unif\{x,w\} & \quad\text{ if $\alpha = b$; } \\
  \unif\{y,w\} & \quad\text{ if $\alpha = c$; } \\
  \unif\{z,w\} & \quad\text{ if $\alpha = d$. }
\end{cases}
\end{equation}
It is clearly stochastically monotone with respect to $(\Ac_0,\Sc_0)$.
We can prove that it is not realizably monotone by contradiction.
Assume that there exists a system $(\XX_\alpha: \alpha\in A_0)$ of
$S_0$-valued random variables which realizes the monotonicity.
Then we have
\begin{align*}
  \Prob(\XX_b = x)
& = \Prob(\XX_a = x, \XX_b = x, \XX_c = w, \XX_d = w)
  = \frac{1}{2}, \\
  \Prob(\XX_c = y)
& = \Prob(\XX_a = y, \XX_b = w, \XX_c = y, \XX_d = w)
  = \frac{1}{2}.
\end{align*}
Therefore, we have $\Prob(\XX_d = w) \ge 1$,
which contradicts the requirement $\Prob(\XX_d = w) = P_d(\{w\}) = 1/2$.
Thus monotonicity equivalence fails for $(\Ac_0,\Sc_0)$.
\qed
\EndExample

In Example~\THn{dia.y.ex}, the dual $\Ac_0^*$ is the diamond again.
By Proposition~\THn{dual}, monotonicity equivalence fails
for both $(\Ac_0^*,\Sc_0)$ and $(\Ac_0,\Sc_0^*)$.
Now let $\Sc$ be any poset of Class~Y.
Since $\Sc$ has either the Y-poset $\Sc_0$ or its dual $\Sc_0^*$ as an
induced subposet,
by Lemma~\THn{s.induce} monotonicity equivalence fails for $(\Ac_0,\Sc)$.
Thus, there exists a system $(\tilde{P}_a,\tilde{P}_b,\tilde{P}_c,\tilde{P}_d)$
of probability measures on $S$ which is stochastically but not realizably
monotone with respect to $(\Ac_0,\Sc)$.
In the next three examples
(Examples~\THn{cycle.y.ex}--\THn{dbow.y.ex}),
we take $\Sc$ to be any poset of Class~Y.

\Example[cycle.y.ex]
Suppose that $\Ac$ has a cycle with height at least $3$.
Then monotonicity equivalence fails for $(\Ac,\Sc)$.

To see this, observe by Lemma~\THn{cycle.height} that $\Ac$ has an
induced cyclic subposet $\Ac' = (a_0,a_1,\ldots,a_{n-1},a_0)$ with
height at least $3$.
Without loss of generality we may assume that
$\Ac'$ has a maximal upward path
$(a_{k'},a_{k'+1},\ldots,a_{n-1},a_0,a_1,\ldots,a_{k})$
with height at least $3$ for some $1 \le k \le k' \le n-1$, as illustrated in
\Equation[cyc.y.fig]
  \Ac' =
  \begin{minipage}{2.5in}\begin{center}
    \includegraphics{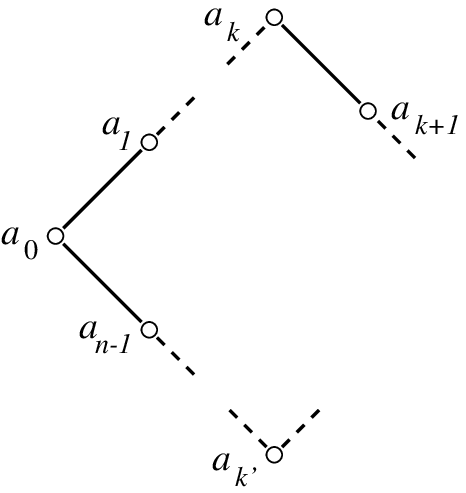}
  \end{center}\end{minipage}
\end{equation}
Note that $k+2 \le k'$, since $a_{k}$ does not cover $a_{k'}$.
Then we can define a system $(P_\alpha: \alpha\in A)$ of probability measures on $S$ by
\begin{equation*}
  P_\alpha := \begin{cases}
    \tilde{P}_b
    & \quad\text{ if $\alpha = a_0$} \\
    \tilde{P}_d
    & \quad\text{ if $a_0 < \alpha$} \\
    \tilde{P}_a
    & \quad\text{ if $\alpha < a_0$} \\
    \tilde{P}_c
    & \quad\text{ otherwise.}
  \end{cases}
\end{equation*}
Then the system is stochastically monotone.

We now show by contradiction that it is not realizably monotone.
Suppose that we have a system $(\XX_\alpha: \alpha\in A)$ of
$S$-valued random variables which realizes the monotonicity.
Then we have $\XX_{a_{k'}} \le \XX_{a_0} \le \XX_{a_{k}}$.
Since $P_{a_{k+1}} = \cdots = P_{a_{k'-1}} = \tilde{P}_c$,
by applying Lemma~\THn{eq.dist.lem} repeatedly
we obtain (almost surely) $\XX_{a_{k+1}} = \cdots = \XX_{a_{k'-1}}$.
Since $a_{k+1} < a_k$ and $a_{k'} < a_{k'-1}$,
we have $\XX_{a_{k'}} \le \XX_{a_{k'-1}} = \XX_{a_{k+1}} \le \XX_{a_{k}}$.
Therefore, the system $(\XX_{a_{k'}},\XX_{a_{0}},\XX_{a_{k+1}},\XX_{a_{k}})$
of $S$-valued random variables realizes the monotonicity
of the system $(\tilde{P}_a,\tilde{P}_b,\tilde{P}_c,\tilde{P}_d)$
in terms of $(\Ac_0,\Sc)$.
But this contradicts Example~\THn{dia.y.ex}.
\qed
\EndExample

\Example[crown.y.ex]
Suppose that $\Ac$ has a $k$-crown $\Ac_k$ as an induced subposet for some $k \ge 3$.
Then monotonicity equivalence fails for $(\Ac,\Sc)$.

To see this, let $\Ac_k$ be as labeled in~\EQn{a.crown},
let $U := \{\alpha\in A: a_0 \le \alpha\}$ be the up-set in $\Ac$
generated by $a_0$,
and let $V := \pideal{b_{k-1}}$ be the down-set in $\Ac$ generated by $b_{k-1}$.
Then we define a system $(P_\alpha: \alpha\in A)$
of probability measures on $S$ by
\begin{equation*}
  P_\alpha := \begin{cases}
    \tilde{P}_b
    & \quad\text{ if $\alpha \in U \cap V$} \\
    \tilde{P}_d
    & \quad\text{ if $\alpha \in U \cap V^c$} \\
    \tilde{P}_a
    & \quad\text{ if $\alpha \in U^c \cap V$} \\
    \tilde{P}_c
    & \quad\text{ otherwise (i.e., $\alpha\not\in U \cup V$).}
  \end{cases}
\end{equation*}
Suppose that $\alpha < \beta$ in $\Ac$.
If $\alpha\in U^c\cap V$, then $P_\alpha = \tilde{P}_a \stocle P_\beta$.
If $\alpha\in U \cap V$, then $\beta\in U$ and $P_\beta$ is either
$\tilde{P}_b$ or $\tilde{P}_d$; thus, $P_\alpha = \tilde{P}_b \stocle P_\beta$.
If $\alpha\not\in U \cup V$, then $\beta\not\in V$
and $P_\beta$ is either $\tilde{P}_c$ or $\tilde{P}_d$;
thus, $P_\alpha = \tilde{P}_c \stocle P_\beta$.
If $\alpha\in U \cap V^c$, then $\beta\in U \cap V^c$,
and $P_\alpha = \tilde{P}_d \stocle \tilde{P}_d = P_\beta$.
In each case that $\alpha < \beta$, we have shown $P_\alpha \stocle P_\beta$.
Therefore, the system is stochastically monotone.

Suppose now that we have a system $(\XX_\alpha: \alpha\in A)$ of
$S$-valued random variables which realizes the monotonicity.
Since $P_{a_0} = P_{b_{k-1}} = \tilde{P}_b$,
by Lemma~\THn{eq.dist.lem} we (almost surely) have $\XX_{a_0} = \XX_{b_{k-1}}$
and therefore $\XX_{a_{k-1}} \le \XX_{a_0} \le \XX_{b_0}$.
Since $P_{a_1} = P_{b_1} = \cdots = P_{a_{k-2}} = P_{b_{k-2}} = \tilde{P}_c$,
by applying Lemma~\THn{eq.dist.lem} repeatedly
we obtain (almost surely) $\XX_{a_1} = \XX_{b_1} = \cdots = \XX_{a_{k-2}} =
\XX_{b_{k-2}}$
and therefore $\XX_{a_{k-1}} \le \XX_{a_1} \le \XX_{b_0}$,
which implies that
the system $(\XX_{a_{k-1}},\XX_{a_0},\XX_{a_1},\XX_{b_0})$
of $S$-valued random variables realizes the monotonicity
of the system $(\tilde{P}_a,\tilde{P}_b,\tilde{P}_c,\tilde{P}_d)$
indexed by the diamond~\EQn{a.diamond}.
But this contradicts the discussion following Example~\THn{dia.y.ex}.
Hence $(P_\alpha: \alpha\in A)$ is not realizably monotone.
\qed
\EndExample

\Example[dbow.y.ex]
Suppose that $\Ac$ has
\Equation[dbow]
  \Ac' =
  \begin{minipage}{1.5in}\begin{center}
    \includegraphics{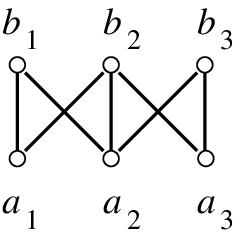}
  \end{center}\end{minipage}
\end{equation}
as an induced subposet.
Then monotonicity equivalence fails for $(\Ac,\Sc)$.

To see this, define a system $\Pc = (P_\alpha: \alpha\in A)$ of probability
measures on $S$ by
\begin{equation*}
P_\alpha := \begin{cases}
  \tilde{P}_b
  & \text{ if $\alpha\in \pideal{b_1}\setminus\pideal{b_3}$ } \\
  \tilde{P}_c
  & \text{ if $\alpha\in \pideal{b_3}\setminus\pideal{b_1}$ } \\
  \tilde{P}_a
  & \text{ if $\alpha\in \pideal{b_1}\cap\pideal{b_3}$ } \\
  \tilde{P}_d
  & \text{ otherwise (i.e., $\alpha\not\in\pideal{b_1,b_3}$). }
\end{cases}
\end{equation*}
Suppose that $\alpha < \beta$ in $\Ac$.
If $\alpha\in\pideal{b_1}\cap\pideal{b_3}$,
then $P_\alpha = \tilde{P}_a \stocle P_\beta$.
If $\alpha\in\pideal{b_1}\setminus\pideal{b_3}$,
then $\beta\not\in\pideal{b_3}$ and $P_\beta$
is either $\tilde{P}_b$ or $\tilde{P}_d$; thus,
$P_\alpha = \tilde{P}_b \stocle P_\beta$.
If $\alpha\in\pideal{b_3}\setminus\pideal{b_1}$,
then $\beta\not\in\pideal{b_1}$ and $P_\beta$
is either $\tilde{P}_c$ or $\tilde{P}_d$;
thus, $P_\alpha = \tilde{P}_c \stocle P_\beta$.
If $\alpha\not\in\pideal{b_1,b_3}$,
then $\beta\not\in\pideal{b_1,b_3}$;
thus, $P_\alpha = \tilde{P}_d \stocle \tilde{P}_d = P_\beta$.
Therefore, the system is stochastically monotone.

Suppose now that there exists a system
$(\XX_\alpha: \alpha\in A)$ of $S$-valued random variables which
realizes the monotonicity.
By Lemma~\THn{eq.dist.lem}, we (almost surely) have
$\XX_{a_1} = \XX_{b_1}$ and $\XX_{a_3} = \XX_{b_3}$.
Thus we have found a system
$(\XX_{a_2},\XX_{a_1},\XX_{a_3},\XX_{b_2})$
of $S$-valued random variables which
realizes the monotonicity of the system
$(\tilde{P}_a,\tilde{P}_b,\tilde{P}_c,\tilde{P}_d)$
indexed by the diamond~\EQn{a.diamond}.
But this again contradicts the discussion following Example~\THn{dia.y.ex}.
Therefore, $(P_\alpha: \alpha\in A)$ is not
realizably monotone.
\qed
\EndExample

\SubSSection[class.y.me]{The proof of Theorem~\THn{class.y.main}}

\SStop
At the beginning of this Section~\SSn{class.y}, we saw
that a non-acyclic poset $\Ac$ can sometimes be enlarged to an acyclic poset $\tilde{\Ac}$.
But such an enlargement is not always possible.
In fact, by Proposition~\THn{non.acyclic},
if $\Ac$ has an induced subposet which is one of
the posets~(i)--(iv) in Proposition~\THn{non.acyclic}
[i.e., (i) the diamond, (ii) a subdivided crown with height at least $3$,
(iii) the $k$-crown for some $k \ge 3$, or (iv) the double-bowtie poset],
then such an enlargement is not possible.
It turns out that a non-acyclic poset can be enlarged to an acyclic poset
if and only if none of the posets~(i)--(iv) in
Proposition~\THn{non.acyclic} is an induced subposet;
this relates the examples in Section~\SSn{class.y.mi} to
Theorem~\THn{class.y.main}.

\Proposition[cyc.enlarge]
Let $\Ac$ be a connected poset.
The following conditions~{\rm (a)}--{\rm (c)} for $\Ac$ are equivalent:
\abcList
\item
there exists an acyclic poset $\tilde{\Ac}$ which has $\Ac$ as
an induced subposet;
\item
any induced cyclic subposet of $\Ac$ is a $2$-crown, and
no induced subposet of $\Ac$ is the double-bowtie~\EQn{dbow};
\item
no induced subposet of $\Ac$ is one of
the posets~{\rm (i)}--{\rm (iv)} in Proposition~{\rm\THn{non.acyclic}}.
\EndList
\EndProposition

\Proof[{\bf Proof of (a) $\Rightarrow$ (c) and (c) $\Rightarrow$ (b).}]
Suppose that there exists an induced subposet $\Bc$ of $\Ac$
which is poset-isomorphic to one of the posets~(i)--(iv)
in Proposition~\THn{non.acyclic}.
If there is an acyclic poset $\tilde{\Ac}$ which has $\Ac$ as an
induced subposet, then $\Bc$ is also an induced subposet of
$\tilde{\Ac}$; by Proposition~\THn{non.acyclic}, this is impossible.
We have thus shown that (a) $\Rightarrow$ (c).

To prove (c) $\Rightarrow$ (b),
observe that a cycle is simply a subdivision of
either the diamond or a crown.
So if $\Ac$ has an induced cyclic subposet $\Bc$ which is not a
$2$-crown, then we can find an induced subposet of $\Bc$
(automatically, of course, an induced subposet of $\Ac$)
that is one of the posets~(i)--(iii) in Proposition~\THn{non.acyclic}.
Thus, the failure to satisfy the condition~(b)
implies the existence of an induced subposet which is
one of the posets~(i)--(iv) in Proposition~\THn{non.acyclic}.
\QED

In preparation for proving (b) $\Rightarrow$ (a),
we introduce a new operation that welds two posets into one, as follows.
Suppose that two posets $\Ac'$ and $\Ac''$ share
a single element $c$ (i.e., that $A'\cap A'' = \{c\}$).
Then the two Hasse diagrams of $\Ac'$ and $\Ac''$
can be drawn in the same plane with their own vertices and arcs
independently except for the vertex $c$ to be shared by the two diagrams;
this introduces a merged diagram on the vertex set $A' \cup A''$.
We call the poset represented by this Hasse diagram
the {\em union of $\Ac'$ and $\Ac''$ joined at $c$\/}
and denote it by $\Ac'\stackrel{c}{\sqcup}\Ac''$.
We list some easily verified properties of the welding operation here:

\numList
\item
If $\Ac'$ and $\Ac''$ are connected, then 
$\Ac'\stackrel{c}{\sqcup}\Ac''$ is  connected.
\item
If $\Ac'$ and $\Ac''$ are acyclic, then 
$\Ac'\stackrel{c}{\sqcup}\Ac''$ is acyclic.
\item
Both $\Ac'$ and $\Ac''$ are induced subposets of
$\Ac'\stackrel{c}{\sqcup}\Ac''$.
\item
Suppose that two posets $\Ac_0'$ and $\Ac_0''$ share a single element
$c$ and that $\Ac'$ and $\Ac''$ (also sharing $c$) are induced subposets of
$\Ac_0'$ and $\Ac_0''$, respectively.
Then $\Ac'\stackrel{c}{\sqcup}\Ac''$ is an induced subposet
of $\Ac_0'\stackrel{c}{\sqcup}\Ac_0''$.
\EndList

We now continue our preparation for the proof of
(b) $\Rightarrow$ (a) in Proposition~\THn{cyc.enlarge}.
Lemma~\THn{lem.split} provides machinery
to split a poset into two smaller ones;
this enables us to devise induction arguments in proving
both Theorem~\THn{class.y.main} and
(b) $\Rightarrow$ (a) in Proposition~\THn{cyc.enlarge}.

\Lemma[lem.split]
Let $\Ac$ be a connected non-acyclic poset.
Suppose that $\Ac$ satisfies the condition~{\rm (b)}
of Proposition~{\rm\THn{cyc.enlarge}}.
Then we can construct a pair $\Ac_0'$ and $\Ac_0''$ of connected posets
(with ground sets $A_0'$ and $A_0''$, respectively) such that,
for some $c$,
\romanList
\item
both $\Ac_0'$ and $\Ac_0''$ satisfy the condition~{\rm (b)}
of Proposition~{\rm\THn{cyc.enlarge}};
\item
$A_0' \cap A_0'' = \{c\}$,
$(A_0' \cup A_0'')\setminus\{c\} = A$,
and $|A_0'|, |A_0''| < |A|$; and
\item
$\Ac$ is the subposet of $\Ac_0 = \Ac_0'\stackrel{c}{\sqcup}\Ac_0''$
induced by $A$.
\EndList
\EndLemma

\Proof
Let $\Gc$ be the collection of all subsets $B$ of $A$
such that the subposet via induced cover subgraph of $\Ac$
on the ground set $B$ is a poset of the form~\EQn{bipart}
for some $m,n \ge 2$.
Note that if a subposet $\Bc$ via induced cover subgraph
is of the form~\EQn{bipart} then $\Bc$ is an induced subposet.
By Lemma~\THn{cons.cycle} and the condition~(b)
of Proposition~\THn{cyc.enlarge},
$\Ac$ has an induced cyclic subposet which is necessarily
a $2$-crown; thus, $\Gc$ is nonempty.
Choose a maximal subset $B_0$ from $\Gc$.
Then let $\Bc_0$ be the subposet of $\Ac$ induced by $B_0$
and label it as in~\EQn{bipart}.

Consider the Hasse diagram of $\Ac$ as represented in the plane.
First remove the arcs from each of $a_1,\ldots,a_m$ to each of $b_1,\ldots,b_n$.
Since the elements $b_1,\ldots,b_n$ are all drawn above
the elements $a_1,\ldots,a_m$,
we can insert a new vertex $c$ above the elements $a_1,\ldots,a_m$
but below the elements $b_1,\ldots,b_n$,
and then add new arcs from $a_1,\ldots,a_m$ to $c$
and from $c$ to $b_1,\ldots,b_n$.
This creates a new poset $\Ac_0$ with ground set $A_0 := A\cup\{c\}$,
as illustrated in
\vspace{0.1in}
\begin{center}
$\Ac_0 = $
\begin{minipage}{2.6in}\begin{center}
  \includegraphics{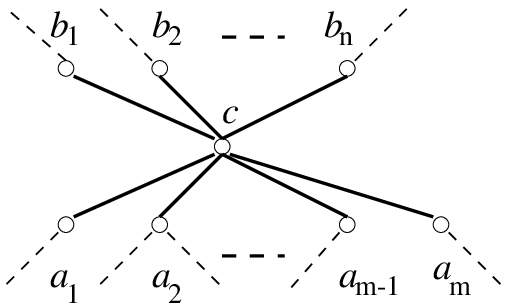}
\end{center}\end{minipage}
\end{center}
The subposet of $\Ac_0$ induced by $A$ introduces
the arc from each of $a_1,\ldots,a_m$ to each of $b_1,\ldots,b_n$,
thus restoring the Hasse diagram of $\Ac$.

We claim that $\Ac_0$ does not have any cycle which contains
an upward path $(a_i,c,b_j)$ with $a_i,b_j\in B_0$.
Granting the claim for the remainder of this paragraph,
we define
\begin{align*}
A_0'  & := \{\alpha\in A_0:
  \text{ a path $(c,a_i,\ldots,\alpha)$ exists in $\Ac_0$
         for some $a_i \in B_0$}\}; \\
A_0'' & := \{\alpha\in A_0:
  \text{ a path $(c,b_j,\ldots,\alpha)$ exists in $\Ac_0$
         for some $b_j \in B_0$}\}.
\end{align*}
By convention, we include $c$ both in $A_0'$ and in $A_0''$.
Let $\Ac_0'$ and $\Ac_0''$ be the subposets of $\Ac_0$ induced by
$A_0'$ and $A_0''$, respectively.
Clearly, $\Ac_0'$ and $\Ac_0''$ are both connected posets.
By observing that $\Ac_0$ is connected, we find $A_0' \cup A_0'' = A_0$.
The claim implies that $A_0' \cap A_0'' = \{c\}$
and that there are no edges between $A' := A_0'\setminus\{c\}$
and $A'' := A_0''\setminus\{c\}$
in the cover graph of $\Ac_0$.
Therefore, $\Ac_0 = \Ac_0'\stackrel{c}{\sqcup}\Ac_0''$,
which implies (iii).
Since $\{a_1,\ldots,a_m\} \subset A_0'$,
$\{b_1,\ldots,b_n\} \subset A_0''$, and $m,n \ge 2$,
we have $|A_0'|, |A_0''| < |A|$, as desired in (ii).
To see (i) for $\Ac_0'$ (the same argument works for $\Ac_0''$),
suppose that $\Ac_0'$ has an induced subposet $\Vc$
which violates the condition~(b) of Proposition~\THn{cyc.enlarge}.
If $V \subset A'$, then let $V' := V$; otherwise,
let $V' := (V\setminus\{c\})\cup\{b_1\}$.
Then the subposet $\Vc'$ of $\Ac$ induced by $V'$ is poset-isomorphic to $\Vc$.
Furthermore, if $\Vc$ is a cycle in $\Ac_0'$, then $\Vc'$ is so in $\Ac$.
But the existence of such an induced subposet of $\Ac$ contradicts
the assumed condition~(b) of Proposition~\THn{cyc.enlarge}.
So (i) holds, and Lemma~\THn{lem.split} is established modulo a proof of
the claim.

Now we show the claim.
To prove this by contradiction, we further enlarge the poset $\Ac_0$ as follows.
We first remove the arcs from $c$ to each of $b_1,\ldots,b_n$
from the Hasse diagram of $\Ac_0$.
Then a new element $c'$ is drawn above $c$ but below
each of $b_1,\ldots,b_n$, and
the arc from $c$ to $c'$ and the arcs from $c'$ to
each of $b_1,\ldots,b_n$ are introduced in the diagram.
This creates a new poset $\Ac_1$, as illustrated in
\vspace{0.1in}
\begin{center}
$\Ac_1 = $
\begin{minipage}{2.6in}\begin{center}
  \includegraphics{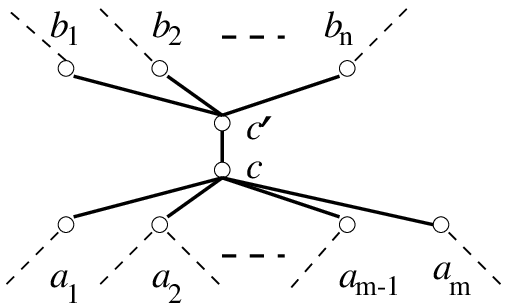}
\end{center}\end{minipage}
\end{center}
Clearly, $\Ac_0$ is an induced subposet of $\Ac_1$.
If we can show that $\Ac_1$ has no cycle which contains the edge
$\{c,c'\}$, then $\Ac_0$ has no cycle which contains an upward path
$(a_i,c,b_j)$ for some $a_i,b_j \in B_0$, establishing the claim and
Lemma~\THn{lem.split}.
Thus, to obtain a contradiction, suppose that $\Ac_1$ has a cycle
which contains the edge $\{c,c'\}$.
By Lemma~\THn{cons.cycle}, we may assume that such a cycle,
say $(a_1,c,c',b_1,u_1,\ldots,u_k,a_1)$, is an induced cyclic subposet of $\Ac_1$.
Then the cycle $(a_1,b_1,u_1,\ldots,u_k,a_1)$ in $\Ac$ is an induced
subposet of $\Ac$, and therefore by condition~(b) of Proposition~\THn{cyc.enlarge}
the induced cyclic subposet $(a_1,b_1,u_1,\ldots,u_k,a_1)$ of $\Ac$
must be a $2$-crown, and therefore $k=2$.
Note that $u_i \not\in B_0$ for $i=1,2$; otherwise,
the cycle $(a_1,c,c',b_1,u_1,u_2,a_1)$ cannot be an induced subposet
of $\Ac_1$.

Write $a_0 := u_1$ and $b_0 := u_2$.
Now consider the comparability in $\Ac$ between $\{a_0,b_0\}$ and $B_0$.
If $a_0$ is comparable with some $a_i$ of $B_0$,
then either $a_0 < a_i < b_1$ or $a_i < a_0 < b_1$,
contradicting our knowledge that $\{a_0,b_1\}$ and
$\{a_i,b_1\}$ are edges in the cover graph of $\Ac$.
Thus $a_0$ is incomparable with each of $a_1,\ldots,a_m$.
Similarly we can see that $b_0$ is incomparable with each of $b_1,\ldots,b_n$.
If $a_0$ is comparable with some $b_j$ of $B_0$ with $j \ge 2$,
then $a_0 < b_j$; otherwise, $b_j < a_0 < b_1$, contradicting the
assumption that $b_1$ and $b_j$ are incomparable.
Suppose that there is an upward path $(v_1,\ldots,v_l)$ in $\Ac$ from
$v_1 = a_0$ to $v_l = b_j$ with $l \ge 3$.
Then it is not hard to see that the cycle
$(a_0,v_2,\ldots,v_l,a_1,b_0,a_0)$ in $\Ac$ is an induced
subposet of $\Ac$ with height $l$,
contradicting condition~(b) of Proposition~\THn{cyc.enlarge}.
Therefore, $b_j$ must cover $a_0$ in $\Ac$.
Note that $a_0$ cannot be comparable with all of $b_2,\ldots,b_n$,
since $B_0\cup\{a_0\} \not\in \Gc$.
Similarly, $b_0$ cannot be comparable with all of $a_2,\ldots,a_m$
(but $b_0$ may cover some of them).
Thus, we can find some elements $a_i$, $b_j$ of $B_0$ so that
$a_0$ is incomparable with $b_j$ and $b_0$ is incomparable with $a_i$.
We have now found the subposet of $\Ac$ induced by $\{a_0,a_1,a_i,b_0,b_1,b_j\}$
to be poset-isomorphic to the poset~\EQn{dbow}.
This contradicts condition~(b) of Proposition~\THn{cyc.enlarge}.
\QED

We now give the proof for (b) $\Rightarrow$ (a) in Proposition~\THn{cyc.enlarge},
by induction over the cardinality of $A$.
The idea of the proof is to build up an acyclic poset
by using Lemma~\THn{lem.split}.

\Proof[{\bf Proof of (b) $\Rightarrow$ (a) in Proposition~\THn{cyc.enlarge}.}]
Suppose that a poset $\Ac$ satisfies the condition~(b).
We will make an inductive argument over the cardinality of $A$.
But if $\Ac$ is acyclic, then the argument is vacuous.
In particular, $\Ac$ with cardinality at most $3$ is acyclic.
Now let $\Ac$ be a connected non-acyclic poset with cardinality $n \ge 4$.
By Lemma~\THn{lem.split}, there exists a pair $\Ac_0'$ and $\Ac_0''$ of
connected posets satisfying (i)--(iii) in Lemma~\THn{lem.split}.
Then, by the induction hypothesis and (i)--(ii) in Lemma~\THn{lem.split},
$\Ac_0'$ and $\Ac_0''$ can be enlarged
to acyclic posets $\tilde{\Ac}'$ and $\tilde{\Ac}''$,
respectively.

Since $A_0'\cap A_0'' = \{c\}$,
the ground sets $\tilde{A}'$ and $\tilde{A}''$ can be given
so that $\tilde{A}'\cap \tilde{A}'' = \{c\}$.
Let $\tilde{\Ac} := \tilde{\Ac}'\stackrel{c}{\sqcup}\tilde{\Ac}''$.
Then $\tilde{\Ac}$ is acyclic.
Furthermore, $\Ac_0 = \Ac_0'\stackrel{c}{\sqcup}\Ac_0''$
is an induced subposet of $\tilde{\Ac}$.
By (iii) in Lemma~\THn{lem.split},
$\Ac$ is an induced subposet of $\tilde{\Ac}$, as desired.
\QED

Now we turn to the proof of Theorem~\THn{class.y.main}.
The proof will parallel that of (b) $\Rightarrow$ (a)
in Proposition~\THn{cyc.enlarge} somewhat.
Let $\Ac$ be a connected poset and let $\Sc$ be a poset of Class~Y.

\Proof[{\bf Proof of Theorem~\THn{class.y.main}.}]
Suppose first that a non-acyclic poset $\Ac$ is
not enlargeable to an acyclic poset.
Then, by Proposition~\THn{cyc.enlarge},
$\Ac$ has an induced subposet $\Bc$ which is
one of the posets~(i)--(iv) in Proposition~\THn{non.acyclic}.
If $\Bc$ is the diamond, then by Lemma~\THn{diamond} $\Ac$ has a cycle
with height at least $3$.
Thus Example~\THn{cycle.y.ex} implies that monotonicity equivalence fails for $(\Ac,\Sc)$.
If $\Bc$ is either the $k$-crown for some $k \ge 3$ or
the double-bowtie poset~\EQn{dbow}, then by
Examples~\THn{crown.y.ex}--\THn{dbow.y.ex} monotonicity equivalence fails for
$(\Ac,\Sc)$.
Suppose now that $\Bc$ is a subdivision of the $k$-crown as
displayed and labeled in~\EQn{a.crown} and has height at least $3$.
Then we may assume that there exists $c_0 \in B$ such that $a_0 < c_0
< b_0$ in $\Bc$.
So we find an upward path $(a_0,\ldots,c',c_0,c'',\ldots,b_0)$ in $\Ac$
from $a_0$ to $b_0$ with height at least $3$.
Since $c_0$ is incomparable in $\Ac$ with each of $a_1,\ldots,a_{k-1}$ and
each of $b_1,\ldots,b_{k-1}$,
no upward path in $\Ac$ from any $a_i$ to any $b_j$ contains either
$\{c',c_0\}$ or $\{c_0,c''\}$ as an edge unless $(i,j) = (0,0)$.
Let $(A,\Ec_{\Ac})$ be the cover graph of $\Ac$.
Then there exists a path $(u_1,\ldots,u_k)$ from $u_1 = c''$ to $u_k = c'$
in the graph $(A,\Ec_{\Ac}\setminus\{\{c',c_0\},\{c_0,c''\}\})$.
Thus, $\Ac$ has a cycle $(c',c_0,c'' = u_1,\ldots,u_k)$ with height at
least $3$.
Therefore, by Example~\THn{cycle.y.ex} monotonicity equivalence fails for $(\Ac,\Sc)$.

Suppose now that a poset $\Ac$ is enlargeable to an acyclic poset.
We will prove that monotonicity equivalence holds
for $(\Ac,\Sc)$ by induction over the cardinality of $A$.
If $\Ac$ is acyclic, then, by Theorem~\THn{a.acyclic},
$\Ac$ is a poset of monotonicity equivalence.
Thus, if $|A| \le 3$, then $\Ac$ is acyclic and therefore a poset of 
monotonicity equivalence.
Now let $\Ac$ be a non-acyclic poset with cardinality $n \ge 4$
and let $(P_\alpha: \alpha\in A)$ be a stochastically monotone system
of probability measures on $S$.
By Lemma~\THn{lem.split}, there exists a pair $\Ac_0'$ and $\Ac_0''$ of
posets satisfying (i)--(iii) in Lemma~\THn{lem.split}.
Let $a_1,\ldots,a_m$ be all the elements covered by $c$,
and let $b_1,\ldots,b_n$ be all the elements covering $c$
in $\Ac_0 = \Ac_0'\stackrel{c}{\sqcup}\Ac_0''$.
Since $\Ac$ is an induced subposet of $\Ac_0$,
we have $P_{a_i} \stocle P_{b_j}$
for all $i=1,\ldots,m$ and all $j=1,\ldots,n$.
By Proposition~\THn{p.insert}, 
we can find a probability measure $P_0$
on $S$ such that $P_{a_i} \stocle P_0 \stocle P_{b_j}$
for all $i=1,\ldots,m$ and all $j=1,\ldots,n$.

Let $P_c := P_0$.
Then we can enlarge the system $(P_\alpha: \alpha\in A)$
to a system $(P_\alpha: \alpha\in A_0)$, maintaining stochastic monotonicity.
Note that the subsystems $(P_\alpha: \alpha\in A_0')$ and $(P_\alpha: \alpha\in A_0'')$
are also stochastically monotone.
Since [by Lemma~\THn{lem.split}(i) and Proposition~\THn{cyc.enlarge}]
$\Ac_0'$ [respectively, $\Ac_0''$] is enlargeable to an acyclic poset,
by the induction hypothesis there is a system $(\XX'_\alpha: \alpha\in A_0')$
[respectively, $(\XX''_\alpha: \alpha\in A_0'')$] of $S$-valued random
variables which realizes the monotonicity of
$(P_\alpha: \alpha\in A_0')$ [respectively, of $(P_\alpha: \alpha\in A_0'')$].
Let $A' := A_0'\setminus\{c\}$ and $A'' := A_0''\setminus\{c\}$.
We can define a probability measure $Q$ on
$S^{A_0} = S^{A'} \times S^{\{c\}} \times S^{A''}$ by
$$
Q(\{(\xx',\xi,\xx'')\})
  := Q'(\xx' | \xi) Q''(\xi, \xx'')
  \quad\text{for $(\xx',\xi,\xx'') \in S^{A'} \times S^{\{c\}} \times S^{A''}$,}
$$
where
\begin{align*}
Q'(\xx' | \xi)
& := P(\XX'_\alpha = \pi_\alpha(\xx')\; \forall\alpha\in A' \:|\:
       \XX_c' = \xi); \\
Q''(\xi, \xx'')
& := P(\XX_c'' = \xi,\: \XX_\alpha'' = \pi_\alpha(\xx'')\; \forall\alpha\in A'').
\end{align*}
Then $Q$ realizes the monotonicity of $(P_\alpha: \alpha\in A_0)$.
By Lemma~\THn{a.induce}, $(P_\alpha: \alpha\in A)$ is realizably
monotone;
thus, monotonicity equivalence holds for $(\Ac,\Sc)$.
\QED

\Remark
In the second part of the proof of Theorem~\THn{class.y.main},
we invoke Proposition~\THn{p.insert}, which requires only $\Sc\not\in\B$.
Thus, we have actually proved that
if $\Sc\not\in\B$ and $\Ac$ is enlargeable to an acyclic poset
then monotonicity equivalence holds for $(\Ac,\Sc)$.
\EndRemark

\SSection[class.z]{Probability measures on a path}

\SStop
In Section~\SSn{class.y} we have seen that, when $\Sc\not\in\B$,
stochastic ordering can be decided from the distribution
function~\EQn{def.d.f}.
In this Section~\SSn{class.z}
we establish that the inverse probability transform~\EQn{def.inv.pr.z}
can be used to realize monotonicity when $\Sc\in\Z$;
this result extends Example~\THn{ex.linear}.
As a result, we will obtain

\Theorem[class.z.main]
Let $\Ac$ be any poset and let $\Sc$ be a poset of Class~Z.
Then monotonicity equivalence holds for $(\Ac,\Sc)$.
\EndTheorem

Let $\Sc$ be a poset of Class~Z.
As we observed in Section~\SSn{me.summary}, $\Sc$
is poset-isomorphic to a path, say $(x_1,\ldots,x_n)$.
So a natural linear order $\le_{n}$ of the path is introduced
by declaring $x_i \le_{n} x_j$ if and only if $i \le j$.
In other words, $(S,\le_{n})$ is a rooted tree with root $x_n$
(see Section~\SSn{prob.on.acyc}).
Note that such a linear order $\le_{n}$ is not consistent in general with
the partial order $\le$ of the poset $\Sc$.
In Figure~\ref{fig.class.z} we give
an example of (a) a poset $\Sc$ of Class~Z and (b) its linear order $\le_{n}$.
For every $x_i \in S$,
the section $\csection{x_i} := \{x_j\in S: x_j \le_n x_i\}$
of the path is either an up-set or a down-set in $\Sc$,
which is obvious pictorially in Figure~\ref{fig.class.z}.
In fact, the linearly ordered set $(S,\le_{n})$ is a rooted tree with
root $x_n$; thus, Lemma~\THn{lem.csec} applies.

For a probability measure $P$ on $S$,
the distribution function $F$ of $P$ is given by~\EQn{def.d.f},
that is, $F(x_i) = P(\csection{x_i})$ for each $x_i \in S$.
Furthermore, we can define
the {\em inverse probability transform\/} $P^{-1}$ from $[0,1)$ to $S$ by
\Equation[def.inv.pr.z]
  P^{-1}(u) := \min\{x_k: u < F(x_k)\}
  \quad\text{ for $u \in [0,1)$,}
\end{equation}
where the minimum is given in terms of the linear order $\le_{n}$.
Then we can state an equivalent condition for
stochastic ordering as the following lemma.

\Lemma[inv.pr.lem]
Let $(P_1,P_2)$ be a pair of probability measures on $S\in\Z$.
Then $P_1 \stocle P_2$ if and only if
\Equation[inv.pr.lem]
  \text{\rm $P_1^{-1}(u) \le P_2^{-1}(u)$ in $\Sc$}
  \quad\text{\rm for all $u \in [0,1)$.}
\end{equation}
\EndLemma

\Proof
Suppose first that $P_1 \stocle P_2$.
Let $F_1$ and $F_2$ denote the distribution functions
of $P_1$ and $P_2$, respectively.
Let $u \in [0,1)$ be fixed,
$x_i := P_1^{-1}(u)$, and $x_j := P_2^{-1}(u)$.
If $x_i = x_j$, then \EQn{inv.pr.lem} obviously holds.
If $x_i <_{n} x_j$, then we have
$$
  F_2(x_k) \le F_2(x_{j-1}) \le u < F_1(x_i) \le F_1(x_k)
$$
for all $x_k$ such that $x_i \le_{n} x_k <_{n} x_j$.
By Lemma~\THn{csec.stoc}, the section $\csection{x_k}$ is a down-set
for every $k = i,i+1,\ldots,j-1$, which implies, by
Lemma~\THn{lem.csec}, that $x_i < x_{i+1} < \cdots < x_{j-1} < x_j$
in $\Sc$.
If $x_j <_{n} x_i$, then we have
$$
  F_1(x_k) \le F_1(x_{i-1}) \le u < F_2(x_j) \le F_2(x_k)
$$
for all $x_k$ such that $x_j \le_{n} x_k <_{n} x_i$.
Again by applying Lemmas~\THn{csec.stoc} and~\THn{lem.csec},
we obtain $x_j > x_{j+1} > \cdots > x_{i-1} > x_i$ in $\Sc$.
In any case, \EQn{inv.pr.lem} holds.

Now suppose that \EQn{inv.pr.lem} holds.
Then we can construct a pair $(\XX_1,\XX_2)$ of $S$-valued random
variables satisfying~\EQn{rv.le}--\EQn{rv.marg} via
$$
  \XX_i := P_i^{-1}(\UU)
  \quad\text{ for each $i=1,2$}
$$
with a single random variable $\UU$ uniformly distributed on $[0,1)$.
By Proposition~\THn{kko}, we have $P_1 \stocle P_2$.
This completes the proof.
\QED

Lemma~\THn{inv.pr.lem} is exactly the property needed
to generalize Example~\THn{ex.linear}
to the case where $\Sc$ is a poset of Class~Z.
We complete the

\Proof[{\bf Proof of Theorem~\THn{class.z.main}.}]
Let $(P_\alpha: \alpha\in A)$ be a stochastically monotone system of
probability measures on $S$.
Let $\UU$ be a random variable uniformly distributed on $[0,1)$.
Then we can construct a system $(\XX_\alpha: \alpha\in A)$ of $S$-valued
random variables satisfying~\EQn{rm.marg} via
$$
  \XX_\alpha := P_\alpha^{-1}(\UU)
  \quad\text{ for $\alpha\in A$.}
$$
By Lemma~\THn{inv.pr.lem}, the system
$(\XX_\alpha: \alpha\in A)$ satisfies~\EQn{rm.mono};
thus, $(P_\alpha: \alpha\in A)$ is realizably monotone.
\QED

{\bf Acknowledgments.}
We thank Keith Crank, Persi Diaconis, Alan Goldman, Leslie Hall,
Robin Pemantle, and Edward Scheinerman for providing helpful comments.

\References

\bibitem{Bogart}  Bogart,~K.~P.  (1996).
{\em Introductory Combinatorics.\/}
Harcourt Brace Jovanovich,
New York.



\bibitem{Fill}  Fill,~J.~A.  (1998).
An interruptible algorithm for perfect sampling via Markov chains.
{\em Ann.\ Appl.\ Probab.\/} {\bf 8} 131--162.


\bibitem{KK} Kamae,~T.\ and Krengel,~U. (1978).
Stochastic partial ordering.
{\em Ann.\ Probab.\/} {\bf 6} 1044--1049.

\bibitem{KKO} Kamae,~T.,\ Krengel,~U.,\ and O'Brien,~G.~L.  (1977).
Stochastic inequalities on partially ordered state spaces.
{\em Ann.\ Probab.\/} {\bf 5} 899--912.


\bibitem{mmthesis}  Machida,~M. (1999).
Stochastic monotonicity and realizable monotonicity.
Ph.D.\ dissertation, The Johns Hopkins University.

\bibitem{invprob}  Machida,~M. (2000).
Inverse probability transform on acyclic posets.
Unpublished manuscript.


\bibitem{Propp-Wilson}  Propp,~J.~G. and Wilson,~D.~B. (1996).
Exact sampling with coupled Markov chains and applications to statistical mechanics.
{\em Random Structures and Algorithms\/} {\bf 9}  223--252.


\bibitem{Ross}  Ross,~D.~A. (1993).
A coherence theorem for ordered families of probability measures on a
partially ordered space.
Unpublished manuscript.

\bibitem{Ruschendorf} R\"{u}schendorf,~L. (1991).
Fr\'{e}chet-bounds and their applications.
In: {\em Advances in Probability Distributions with Given
Marginals.\/}
Eds. Dall'Aglio,~G. Kotz,~S. and Salinetti,~G.
151--187. Kluwer Academic Publishers.

\bibitem{Shortt} Shortt,~R.~M (1989).
Combinatorial methods in the study of marginal problems over separable
spaces.
{\em J.\ Math.\ Anal.\/} {\bf 97} 462--479.


\bibitem{Stanley} Stanley,~R.~P.  (1986).
{\em Enumerative Combinatorics\/}, Volume~1.
Wadsworth \& Brooks/Cole,
Monterey, California.

\bibitem{Strassen} Strassen,~V. (1965).
The existence of probability measures with given marginals.
{\em Ann.\ Math. Statist.\/} {\bf 36} 423--439.

\bibitem{Trotter} Trotter,~W.~T.  (1992).
{\em Combinatorics and Partially Ordered Sets: Dimension Theory.\/}
The Johns Hopkins University Press,
Baltimore, Maryland.

\bibitem{Vorobev}  Vorob'ev,~N.~N. (1962).
Consistent families of measures and their extensions.
{\em Theory\ of\ Probab.\ Appl.\/} {\bf 7}  147--163.

\bibitem{West}  West,~D.~B.  (1996).
{\em Introduction to Graph Theory.\/}
Prentice Hall,
New Jersey.


\EndReferences

\par\noindent
\textsc{James Allen Fill \\
Department of Mathematical Sciences \\
The Johns Hopkins University \\
Baltimore, MD 21218-2682 \\}
\texttt{jimfill@jhu.edu}

\vspace{0.25in}

\par\noindent
\textsc{Motoya Machida \\
Department of Mathematical Sciences \\
The Johns Hopkins University \\
Baltimore, MD 21218-2682 \\}
\texttt{machida@mts.jhu.edu}

\end{document}